\documentclass{slxx}

\usepackage{slsec} 
\usepackage{stud_log}
 \usepackage{xxslfoot} 
 \usepackage{slthm}

\usepackage{
amsfonts,
amsmath,
latexsym,
amssymb,
enumerate,
verbatim,
mathrsfs,
graphicx,
centernot,
accents, 
color
}

\usepackage[all]{xy}

\usepackage{tikz} \usetikzlibrary{arrows}\usepackage{graphicx}

\newcommand{\labbel}[1]{\label{#1} [[{\bf #1}]]}  
\renewcommand{\labbel}{\label}

 \definecolor{reed}{RGB}{0,0,25} 

\newcommand{\arxiv}[1]{{\color{reed}#1}}

\newcommand{\K}{K} 
\newcommand{\Ku}{K}
\newcommand{\CC}{C}

\newcommand{\vdashcup}{\mathrel{\vdash\hspace{-5 pt}^{^{_\cup}}}}

\newcommand{\sharr}{\text{\tiny{$\rightarrow$}}}
\newcommand{\larr}{\text{\tiny{$\leftarrow$}}}

\newcommand{\nup}{{\centernot\uparrow}}
\newcommand{\da}{{\downarrow}}

\newtheorem{theorem}{Theorem}[section]
\newtheorem{lemma}[theorem]{Lemma}

\newtheorem{proposition}[theorem]{Proposition} 
 
\newtheorem{corollary}[theorem]{Corollary}

\newtheorem*{claim*}{Claim}

\newtheorem*{theorem*}{Theorem}
\newtheorem*{proposition*}{Proposition}
\newtheorem*{corollary*}{Corollary}
\newtheorem*{lemma*}{Lemma}
\newtheorem*{scholion*}{Scholion}

\theoremstyle{definition}
\newtheorem{definition}[theorem]{Definition}

\newtheorem{problem}[theorem]{Problem} 
\newtheorem{problems}[theorem]{Problems}

\theoremstyle{remark}
\newtheorem{remark}[theorem]{Remark}
\newtheorem{remarks}[theorem]{Remarks}
\newtheorem*{remark*}{Remark}
\newtheorem*{remarks*}{Remarks}
 
\newtheorem{example}[theorem]{Example}

\newtheorem*{observation*}{Observation}

\newtheorem{remarkx}{Remark}[theorem]
\newtheorem{observationx}[remarkx]{Observation}

\theoremstyle{definition}
\newtheorem{problemx}[remarkx]{Problem}

\numberwithin{equation}{section}

\HeadingsInfo{Paolo Lipparini}
{A model theory of topology}

\begin{document}

 \AuthorTitle{Paolo Lipparini}
{A model theory of topology}

\begin{abstract}
 An algebraization of the notion of topology
has been proposed more than seventy years ago
 in a classical paper by  McKinsey and Tarski, leading
to an area of research still active today,  with connections
to algebra, geometry, logic and 
many applications, in particular, to modal logics.
In McKinsey and Tarski's setting
 the model theoretical  notion of homomorphism does not correspond
to the notion of continuity.
We notice that the two notions correspond if instead
we consider a preorder relation $ \sqsubseteq $ 
defined by $a \sqsubseteq b$ if
 $a$ is contained in the topological closure of $b$,
for $a,b$ subsets of some topological space.
 
A \emph{specialization poset}
is a partially ordered set endowed with a further
coarser preorder relation $ \sqsubseteq $.  We show that 
every specialization poset can be embedded in
the specialization poset naturally associated to some topological space,
where the order relation corresponds to set-theoretical inclusion.
Specialization semilattices are defined in an analogous way
and the corresponding  embedding theorem  is proved.
Specialization semilattices have
 the amalgamation
property. Some basic topological facts and notions are recovered in this
apparently very weak setting. 
The interest of these structures 
arises from the fact that they also occur in many rather disparate contexts,
even far  removed from topology. 
\end{abstract} 


\arxiv{
\bigskip 

The present ArXiv version contains a few additional
comments which are not present in the published version
(\emph{A Model Theory of Topology}, Studia Logica,
Volume 113, pages 225--259, 2025). Some further
bibliographical references have also been added.

Further results about the notions studied here appear 
in \cite{mttlib,mttna,mttcong}.
Quite surprisingly, a somewhat simpler
and neater theory arises if a ternary
relation corresponding to $z \subseteq y \cup K w$
is taken into account \cite{rsaif}.

In \cite{mttmult}
$n+1$-ary relations are studied, 
whose intended interpretation is
$z \subseteq Kw_1 \cup Kw_2 \cup\dots\cup Kw_n $, 
with possible connections with multi-posets \cite{rump}.  
We have not yet studied what happens merging the
ideas of \cite{rsaif}
and \cite{mttmult}.

In a slightly different direction, further topological notions preserved
under continuous images in the sense studied here
are discussed in \cite{cs,hyp,csb,cp}.  
In particular, see \cite{hyp} for further motivations and
connections with region-based theory of space, event 
 structures in the theory of concurrent 
systems in computer science, representations of median  
graphs, the theory of intersection graphs and hypergraphs,
network science. The manuscript  \cite{rsaif} tries to merge and unify the
different  approaches.
}

\section{Introduction} \labbel{intro}

 As is well known,  a topology  can be equivalently characterized 
by means of the corresponding Kuratowski closure operator $ \Ku $, e.g.\  
\cite[Ch.\ II, Sect.\ 3]{Nag}. 
The characterization naturally lends itself to an algebraization
for the notion of topology.
In the seminal paper 
 ``The algebra of Topology'' \cite{MT} McKinsey and Tarski introduced 
closure algebras, which are Boolean algebras endowed with
an additional  operation $ \K $  satisfying the formal properties of 
topological closure for subsets of a topological 
space\footnote{See also Sikorski \cite[footnote 1]{Si} and
footnote 1 in Kuratowski    
\cite[Ch.\ I, Sect.\ 4.I, p.\ 20]{K} for reference to earlier works
on the subject.
The area of research is still very active today,  with connections
to algebra \cite{Es}, geometry \cite{vBB,E}  and logic \cite{RS},
 particularly,  modal logics \cite{Ar,vBB,GKWZ,GM}.
Recent works on the subject include, just to mention some,
variations on McKinsey and Tarski Theorem   \cite{BBLM,La},
sometimes involving deep set theoretical combinatorial problems \cite{LV},
new duality theories \cite{BMM,Cel},
extensions to more powerful languages \cite{CPZ,GH}
and an accurate study of existentially closed closure
algebras, possibly with the additivity condition removed \cite{Sc}. 
 See the quoted sources 
for more details and references.}. 

In the present paper we aim to approach some 
foundational aspects of 
McKinsey and Tarski's method.
The motivation comes from the observation that
the construction from \cite{MT} is not 
``functorial''.
Algebraically, a homomorphism $\psi$
of closure algebras satisfies
$\psi( \K z)= \K  \psi(z)$.
On the other hand, if
$X$, $Y$ are topological spaces, 
  $ \varphi $ is a function from 
 $X$ 
to  $Y$
and $ \varphi ^\sharr $   is the corresponding image function
from the power set $\mathcal P(X)$ to 
$\mathcal P(Y)$,
then    
$ \varphi $ is continuous if and only if  
$ \varphi ^\sharr ( \K z) \subseteq  \K  \varphi ^\sharr(z)$,
for every $z \subseteq X$, where $  \K $ is topological closure.
The converse inclusion 
holds if and only if  $ \varphi $ is a closed map
\cite[Exercise 1.4.C]{En}.
Not every continuous map 
is closed.

The correspondence 
between topological spaces and closure algebras is not even contravariant,
with respect to  the algebraic notion of homomorphism.
If $X, Y$ are topological spaces, $ \varphi : X \to Y$ is
 a function and $ \varphi ^\larr $ denotes
 the preimage function, 
then $ \varphi $ is continuous if and only if 
 $  \K  \varphi ^\larr (z) \subseteq \varphi ^\larr ( \K  z)$,
for every $z \subseteq Y$ 
\cite[Sect.\ 13.IV(2)]{K}, but equality holds if and only if $\varphi$  is also open
\cite[Exercise 1.4.C]{En}. Open and continuous
functions are sometimes called \emph{interior maps} \cite{RS}. 
As above, not every continuous function is open. 
Thus, in order to exploit the topological content of 
closure algebras, an appropriate notion of a morphism
is  a Boolean algebra homomorphism satisfying 
$ \K  \psi(z) \leq \psi( \K z) $. Such morphisms have been
variously called  \emph{continuous morphisms}, 
\emph{semi-homomorphisms} 
\cite{BMM,Si} 
and possibly even  in different ways.
With such a choice of morphisms, functoriality
is restored  \cite{BMM}.
The approach turns out
to be very useful; on the other hand,
 it looks quite unnatural from the algebraic point of view.

We  observe that if, instead, we consider a
binary relation $ z \sqsubseteq w $,
to be interpreted as
\begin{equation}\labbel{1}      
\text{$ z \sqsubseteq w $ \quad if
\quad  $z$ is contained in the closure of $w$,  }
 \end{equation}
for $z, w \subseteq X$, then model-theoretical homomorphisms
correspond exactly  and covariantly  to continuous functions.
The above observation is an immediate generalization of the well-known fact
that a function $\varphi$   is continuous between 
topological spaces if and only if 
$\varphi$   preserves the adherence relation
$x \in  \K z$ defined between a point $x $ and a subset
$z $ \cite[Sect.\ 13]{K}.

Thus we are led to consider 
 algebraic and relational structures which are associated to topological spaces 
and are preserved by  image functions
associated to continuous maps\footnote{We mean at the basic level
of the fundamental definitions,
not in the sense of algebraic topology, homology, etc.}.
Since image functions preserve unions
but not necessarily intersections or complements, 
we study join-semilattices
with a coarser preorder $ \sqsubseteq $
satisfying the condition
\begin{align}
\labbel{s3intro}    \tag{S3}
 &a \sqsubseteq b  \ \&\  a_1 \sqsubseteq b
 \Rightarrow 
 a \vee a_1 \sqsubseteq b, 
\end{align}
which is obviously satisfied in our intended example,
where $ \sqsubseteq $ is given by \eqref{1}.
The preorder $ \sqsubseteq $ 
 shall be called a \emph{specialization} and 
 the above  structures  \emph{specialization semilattices}.
The name \emph{specialization} comes from the fact that 
this is a generalization of the \emph{specialization order}
between points of a topological space. See item 1.\ in Section \ref{exasub}
 below for details.  

Our main technical result Theorem \ref{sptop}
  asserts that  every specialization semilattice 
can be embedded into (the structure
associated to) some topological space.
In particular,  our axiomatization captures exactly 
the universal theory valid for subsets of a  topological space 
in the language $\{ \vee, \sqsubseteq  \}$. 
If we consider only those properties holding 
in the language $\{ \leq, \sqsubseteq  \}$,
where $\leq$ is interpreted as inclusion in the motivating example
of topological spaces, we will speak of specialization posets. 
In detail, a \emph{specialization poset}
is a partially ordered set together with 
a coarser preorder. A corresponding representation theorem
for specialization posets is proved in Proposition \ref{embp}.

Summarizing, we have been
searching for a theory which could speak of topological spaces 
and which furthermore satisfies some quite stringent requisites.
  \begin{enumerate}[(R1)]
    \item   
 The theory is expressible as a first-order theory,
preferably using universal sentences. In particular,
we do not deal with infinitary operations, or second-order properties.
\item
 The original notion of
morphism for topological spaces, that is, continuous map,
should correspond to the notion of homomorphism in the
algebraic or  model-theoretical sense,
 and:
\item
 The correspondence between such morphisms should be 
covariant, not contravariant, with respect to the direction
considered in the original setting.
 \end{enumerate}

As well-known, just dropping
any one of requirements (R1) - (R3),
 many successful and consolidated
approaches are viable\footnote{The most classical approach is 
Stone duality \cite{Sto}, leading 
(among others \cite{J})
to various point-free
generalizations of topology \cite{PP,PP2,V}.
McKinsey and Tarski's original
 approach \cite{MT}, together with further developments  \cite{JT,MT2,RS},
 has been shown to provide a 
simple decidable formalism to talk about topology,
which is frequently called ``the modal logic of 
topology'' \cite{Ar}.
Other approaches include
topological model theory 
\cite{BF,CKc,Z},
different axiomatizations \cite{Ed,Se,Pa} and
the study of further  mathematical structures
related to topology \cite{DC,E,MP},
possibly from a constructive point of view \cite{SG}.}.

\arxiv{Thus our project is mainly motivated by  aesthetical taste,
possibly influenced by the model-theoretical background of the author.
Namely, we do not aim  to  create
still another alternative theory; rather, we just  ask 
what happens if we pursue the  model theoretical point
of view in its entirety (and in the first-order setting).
 }

While originally we have pursued only the above-described minimal objective, 
namely, to detect which structural parts of topological spaces
are preserved under model-theoretical homomorphisms
in the first-order setting,
we subsequently realized that the resulting structures  
turn out to be of independent interest.
Though the idea has come out of mere curiosity,
 we ended up with theories which, at least in 
the author's opinion, have some quite remarkable features.
  \begin{enumerate}[(F1)]
    \item  
First of all, the  theories have arisen independently
in many different and distant settings. See
Sections \ref{motiv} and  \ref{exasub}  for a list of examples.
The fact that many applications exist in very different fields,
sometimes far removed from topology,
make the theories  interesting for themselves, independently
of the topological interpretation.

\item
The theories are simply described and axiomatized;
representation theorems are easily proved 
(Theorems \ref{propcomp}, \ref{sptop},   
Proposition \ref{embp} and \cite{mttlib}) 
and do not need the axiom of choice.
In particular, specialization semilattices admit a quite clean structure theory, 
 have the amalgamation property, a Fra\"\i ss\'e model
 and a model-completion. In passing, note that,
on the other hand, the theory of closure algebras 
has no   model-companion \cite{ecca,Sc}.

\item  
Some  topological facts and notions  can be reconstructed
even in this apparently very weak setting, for example,
 some properties
of compactness are retrieved.  See Section \ref{recov}.
\arxiv{As another  example, closures of elements need not necessarily exist
in a specialization semilattice; however, once some element has a closure, this closure
can be proved to satisfy many of the usual properties:
 see Lemma \ref{lemlemlem}.  }

\arxiv{\item  The theories help in clarifying the connections between ``full'' closure
and specialization, either in the form of ``generalized adherence'',
or  in the sense of ``contained in the hull''. In particular,
there are canonical and quite easily described embeddings
of a specialization semilattice into a closure semilattice,
Theorem \ref{propcomp} and \cite{mttlib,mttna}. 
On the other hand, we describe many situations
in which specialization appears a more natural notion than
closure: see  Section \ref{exasub} and 
a few additional comments in the expanded introduction
of \cite{mttlib} 
in arXiv:2201.09083v2.}
 \end{enumerate} 

\arxiv{For short, while the motivations are from logic and foundations,
it seems that the resulting theories have mainly an algebraic and
order-theoretical interest, with many applications, sometimes even outside
mathematics.}
 
The paper is divided as follows. 
In Section \ref{motiv}
we describe the motivating example in more detail.
\arxiv{To every topological space we associate a specialization
semilattice and we check that continuous functions
between topological spaces correspond exactly
to homomorphisms between specialization semilattices.  
The correspondence works more  generally
for closure spaces, actually, for partially ordered sets
with a closure operation.  }
In Section \ref{specsec} we present the 
actual definitions of specialization semilattices
and posets and give a few elementary consequences.
Section \ref{exasub} presents more examples
of specialization semilattices and posets, sometimes
appearing in scientific contexts quite away from topology.
\arxiv{The section is not necessary in order to understand
the remaining parts of the paper. }
In Section \ref{embsec} we prove our main result 
Theorem \ref{sptop}, to the effect that every specialization semilattice 
can be embedded into the specialization semilattice associated 
to some topological space. 
\arxiv{The proof is divided in various steps.
First in Subsection \ref{prin} we study principal specialization semilattices,
roughly, those specialization semilattices in which
a notion of closure can be 
defined. In Subsection \ref{emb}   
 we prove that every specialization semilattice can be embedded into a 
principal specialization semilattice, then in Subsection \ref{univcs}
we complete the proof of Theorem \ref{sptop}.
Namely, we show that 
every principal
specialization semilattice can be embedded into the specialization semilattice 
associated to some topological space. 
Thus the universal theory
of specialization semilattices is the universal theory    
of topological spaces in the language of specialization semilattices. }
In Subsection \ref{sp} we prove the corresponding result
for specialization posets, while in Subsection \ref{ap}
we show that the theory of specialization semilattices has
the amalgamation property, hence Fra\"\i ss\'e limits and model-completion. 
In Section \ref{recov} we show that some topological notions
can be expressed in the language of specialization posets, and
retrieve some properties of compactness.  Section \ref{fr} contains
 further remarks and
problems.
\arxiv{In Section \ref{me} we present more examples and
some counterexamples. }

\arxiv{
\subsection*{Prerequisites.}  
While our original motivation is mainly model-theoretical,
we do not assume a specific knowledge of model theory
from the reader and (apart from Subsection \ref{ap})
 only scattered results here mention
or use model-theoretical notions.
The model-theoretical notions of homomorphism
and embedding are an exception, but they are fully
explained in detail in Subsection \ref{homomorphisms}.
We have tried to make the paper as self-contained as possible.
Only a minimal mathematical background is assumed,
essentially, some basic notions of topology, order-theory and algebra.
A few comments are expressed in logical terminology,
but they are not necessary to understand the remaining parts of the paper.
When dealing with logic, we always work in the setting of first-order 
``classical'' logic, i.e., finitary, two-valued, with
 only Boolean connectives, accepting the 
law of excluded middle\dots }

\section{The motivating example} \labbel{motiv} 

 The topology $\tau$ of some topological space $(X, \tau )$
can be equivalently described
by specifying the family of its closed 
subsets, equivalently, as recalled in the introduction,
its associated 
closure operation $ \Ku _ \tau $.
The operation $ \Ku _ \tau $  sends
a subset $a$ of $ X$ to the intersection of all the closed
subsets of $X$ which contain $a$. 
Clearly, $ \Ku _ \tau $ is determined by
$\tau$. Conversely, given a function
$ \Ku : \mathcal P(X) \to \mathcal P(X)$
satisfying $ \Ku   \emptyset = \emptyset $
(\emph{$ \Ku  $ preserves $ \emptyset $}), $ \Ku  x \supseteq x$  
(\emph{inflationary}, or \emph{extensive}), 
$ \Ku   \Ku  x =  \Ku  x$   (\emph{idempotent}) and 
$ \Ku  (x \cup y) =  \Ku  x \cup  \Ku  y$ (\emph{additive}, or 
\emph{topological}),
the family $\{ a \subseteq X \mid  \Ku  a=a \}$  
is the family of closed sets for some topology. 
The above constructions are one the inverse
of the other \cite[Ch.\ II, Sect.\ 3]{Nag}, hence they provide equivalent descriptions
for a topology.
When no risk of confusion is possible, we shall simply write
$X$ in place of $(X, \tau )$ and 
$ \Ku  $ in place of $ \Ku  _ \tau $.

Essentially all the arguments in the present note
work in the more general context of closure spaces,
which are like topological spaces, except
that the assumptions about the family of closed sets
are relaxed. 
\arxiv{As we will mention,
closure spaces have applications in many distant mathematical fields. }

\begin{definition} \labbel{closp}   
A \emph{closure space}
is a set $X$ together with a family $\mathcal F$ 
of subsets of $X$, 
such that that $\mathcal F$ is closed under
arbitrary intersections.
We assume that $X \in \mathcal F$ 
(this is redundant if one assumes that $X$ is the intersection
of the empty family).
Members of $\mathcal F$  are also called \emph{closed (sets)}.
Thus we leave out the assumption
that $\mathcal F$ is preserved  under finite unions, an assumption
holding in topological spaces. We are not asking that 
$\emptyset $ is closed, either.

 As in the  case
of topologies, closure spaces
can be equivalently characterized  by the 
associated \emph{closure operator} $ \K $.
In the case of closure spaces $ \K $ 
is assumed to be only
inflationary, idempotent and \emph{isotone}; this last 
condition means that 
$x \subseteq y$ implies $ \K x \subseteq  \K y$,
equivalently, $ \K (x \cup y) \supseteq   \K x \cup  \K y$.
A closure space defined in terms of a closure operator
is a topological space if and only if the operator 
is additive and preserves $ \emptyset $. 
When no risk of confusion
might arise, we will simply say that $X$ is a closure space,
with no mention of $\mathcal F$ or $ \K $.
\arxiv{The terminology about closure spaces
is not uniform in the literature
and the name of the notion itself
greatly varies according to the author or to the field of research.
See \cite{E} for further details, a historical survey, 
further references and, in particular, 
 \cite[p. 163]{E} for a picture. }

If $X$ and $Y$ are closure spaces, then, exactly as in the case
of topological spaces, a function $\varphi: X \to Y$
is \emph{continuous} if and only if 
$ \varphi ^\sharr ( \K a)  \subseteq  \K \varphi^\sharr(a) $,
for every $a \subseteq X$.
Equivalently, $\varphi$  is continuous if and only if  
the preimage
of each closed subset of $Y$ is a closed subset of $X$.
The equivalence is proved exactly as in the case of 
topological spaces,  
 e.g.\  \cite[Sect.\ 13.IV]{K} or 
\cite[Proposition 1.4.1]{En}.
\arxiv{For the
reader's convenience,
 we present explicit details.}
 \end{definition}

\arxiv{\begin{lemma*} \labbel{contt}
If $X$ and $Y$ are topological spaces, or just closure spaces,
and  $\varphi $ is  function from $  X  $ to $Y$,
then $\varphi$  is continuous if and only if 
the preimage
of each closed subset of $Y$ is a closed subset of $X$.
 \end{lemma*} 

\begin{proof}
To prove the result in the general case
of closure spaces, suppose that $\varphi$ is continuous
and  $c$ is
closed in  $Y$.
If, by contradiction, the preimage $ \varphi  ^\larr (c)$ is not closed in 
$X$,
there exists  $x \in  \K  \varphi  ^\larr (c)$,
such that $ x \notin \varphi  ^\larr (c)$. 
Since $\varphi ^ \sharr (  \K  \varphi  ^\larr (c)) \subseteq 
    \K  \varphi ^ \sharr( \varphi  ^\larr (c)) \subseteq  \K  (c)=c $, then
$\varphi(x) \in c$, a contradiction.

In the other direction, if preimages 
of closed are closed, then  $ \varphi ^\larr( \K  \varphi^ \sharr (b)) $
is a closed  containing 
$  \varphi ^\larr( \varphi^ \sharr (b)) \supseteq b $,
hence $ \varphi ^\larr( \K  \varphi^ \sharr (b)) $ contains $ \K (b)$.
Thus
$ \varphi ^ \sharr (\varphi ^\larr( \K  \varphi^ \sharr (b)))$
contains $  \varphi^ \sharr ( \K b)$
and the conclusion follows from the fact that
$ \varphi ^ \sharr (\varphi ^\larr( \K  \varphi^ \sharr (b)))
\subseteq   \K  \varphi^ \sharr (b)$. 
 \end{proof}    

Continuity between 
closure spaces seems to have received less attention than it deserves,
since there are many significant examples. }
Non-topological closure
 spaces and continuity 
between them naturally arise
also in purely topological contexts, see \cite{C}
for examples. 

Closure spaces have found many applications in very disparate
settings, with varied and occasionally clashing terminology.
Actually, it is quite difficult to collect all instances of applications
of closure spaces.
An ample discussion is presented
in \cite{E},  with illuminating figures and a detailed historical background,
 highlighting applications, among others,
 to ordered sets, lattice theory,
logic, algebra,  topology, quantum physics
and connections with category theory.
The notion of a closure operation
is formally the same \cite{Jan,W} 
as the notion of an \emph{abstract} (not necessarily finitary)
consequence operation.  See item 4.\ in 
Section \ref{exasub} below for more details.
 A useful set of references to applications
in computer science, notably, in the semantic area, 
can be found in \cite{R}. 
A
 comprehensive survey of applications to
conceptual data analysis, relational database theory
and other areas can be found in \cite{CLM,CM}.
For applications to universal algebra, the reader might
consult \cite{B}. 

Model theoretical properties
of the algebraic analogue of closure spaces,
that is, Boolean algebras with a (not necessarily additive)
closure operation are studied in \cite[Section 8]{Sc}.
Earlier results have been obtained by the Italian school, e.g.\ 
\cite{Se}. 
Just like topological spaces and closure algebras
furnish an algebraization for the modal system S4,
closure spaces and their algebraization are the counterpart
of a non-normal ``monotonic''  modal system. 
See  \cite{Fr,Gr,MC}
for further information about monotone modal logics
and \cite[Footnote 8]{Sc} for the connections 
with Boolean algebras with a nonadditive
closure operation.  
Other useful references 
about closure spaces
are, among many others, \cite{Bly,GHKLMS}. 

\begin{remark} \labbel{anchclos}    
In connection with closure spaces, let us mention that
many remarks from the introduction 
concerning continuous functions and associated homomorphisms
apply to closure spaces, as well.
A typical example of a closure space 
arising in an algebraic context is the closure space
of normal subgroups of a group $\mathbf G$.
Here the underlying set is $G$ and the closed subsets
are (the domains of) normal subgroups.
If $\varphi:\mathbf G \to \mathbf H$    
is a group homomorphism, then 
the  image function 
$\varphi ^\sharr:
\mathcal P(G) \to \mathcal P(H)$ is continuous
between the associated closure spaces, but not necessarily
a homomorphism of closure spaces. In fact, the image of
a normal subgroup is not necessarily normal.
The same remark applies to a general algebraic structure
$\mathbf A$ and the associated closure space on $A \times A$
of congruences on $\mathbf A$ \cite{B}.
A similar remark applies also to deduction systems,
see item 4.\ in Section \ref{exasub} below.

On the other hand, if we associate to a group 
(more generally, to an arbitrary algebraic structure) the closure space of
subgroups (subalgebras), then a homomorphism
   $\varphi$  does induce a homomorphism
$\varphi^\sharr$ of closure spaces.  
\end{remark}

As usual, \emph{poset}
is an abbreviation for \emph{partially ordered set}. 
\emph{Semilattices}   here shall be always  intended as \emph{join-semilattices},
namely, to a semilattice $(S, \vee)$
there is associated the partial order $\leq$ on $S$
defined by 
$a \leq b$ if and only if $a \vee b =b$.
When we mention a partial order in reference to some 
semilattice, we shall always mean the order $\leq$ 
introduced above. 
The symbols 
$ \vee$ and  $\wedge$ shall always be used 
to denote joins and meets in posets, semilattices or lattices. 
They should not be confused with the logical propositional
operators of conjunction and (inclusive) disjunction, which shall be denoted by
 ``$\&$''  and  ``or''.

\begin{definition} \labbel{asso}   
Let $X$
be a topological space, or just a closure space.
The \emph{specialization semilattice
$\mathbf S(X)$ associated to $X$}
is the model 
$( \mathcal P(X),   {\cup}, {\sqsubseteq} )$,
where $ \sqsubseteq $ 
is the binary relation on 
$\mathcal P(X)$
defined by
$ a \sqsubseteq b$
if  $a \subseteq  \K b$.
 Here $a$ and $b$ vary among subsets of $X$  
and $ \K $ is closure.
The \emph{specialization poset
$\mathbf P(X)$ associated to $X$}
is  $( \mathcal P(X), {\subseteq} , {\sqsubseteq} )$,
where $ \sqsubseteq $ is defined as above.

\arxiv{When it is necessary to make explicit reference to the topology
$\tau$ on $X$, we will write 
$\mathbf S(X, \tau )$ and $\mathbf P(X, \tau )$.
In the case of a closure space we shall consider $\tau$ 
as the family of the complements of members of $\mathcal F$. }
 
Explicit  definitions for
 abstract notions of  specialization  semilattices and  posets
will be given in Definition \ref{spsem} below, 
by means of a few natural conditions.
We will show in Theorem \ref{sptop} and Proposition \ref{embp}
 that any structure satisfying 
the conditions in Definition \ref{spsem} is isomorphic
to a substructure of  $\mathbf S(X)$ or
of $\mathbf P(X)$, as introduced above, for some topological space $X$.

Our main interest in the above notions originates from
the next  proposition, which  is essentially known,
see e.g.\ the proof of \cite[Theorem 2.4]{Sla}.
Theorem 2.4 in \cite{Sla} is stated for complete semilattices,
but the semilattice structure is not involved
in those parts of the proof related to the next proposition.
In any case, we will present full details below for the
reader's convenience. 

Here homomorphisms  are always intended in the 
model-\hspace{0 pt}theoretical sense \cite{H}.
\arxiv{See
Subsection \ref{homomorphisms}  below for explicit details 
in the special cases at hand.  }
Recall that a continuous map 
$\iota: X \to Y$ between topological spaces 
is a \emph{(topological) embedding}
 if $\iota$ induces an homeomorphism
from $X$ to $\iota ^\sharr (X)$,
considered as a subspace of $Y$.    
In particular, an embedding is an injective function.
If $Y$ is just a closure space and $Z \subseteq Y$,
then, exactly as for topological spaces,  $Z$ inherits the structure 
of a closure space by taking as closed subsets of $Z$ 
the subsets of the form $Z \cap C$, with $C$ closed in $X$. 
Embeddings of closure spaces are defined as above.
\end{definition}

\begin{proposition} \labbel{propcont}
Suppose that $X$ and $Y$ are topological spaces,
or just closure spaces, and
$\varphi: X \to Y$ is a function.
Then the following conditions are equivalent.
 \begin{enumerate}[(1)] 
  \item  
 $\varphi$  
is continuous from $X$ to $Y$;
\item   
the image function $\varphi ^\sharr 
: \mathcal P(X) \to \mathcal P(Y) $
is a homomorphism  of specialization semilattices
from   $\mathbf S(X)$ to $\mathbf S(Y)$;
\item   
the image function
$\varphi ^\sharr  $
is a homomorphism  of specialization posets
from   $\mathbf P(X)$ to $\mathbf P(Y)$.
   \end{enumerate}

The equivalences still hold if we replace everywhere
``continuous'' and  ``homomorphism'' with ``embedding''.
\end{proposition}

  \begin{proof}
The image 
function $\varphi ^\sharr $
satisfies 
$\varphi ^\sharr (a \cup b) = \varphi ^\sharr (a) \cup \varphi ^\sharr (b)$,
for every $a, b \subseteq X$, with no further special assumption,
hence $\varphi ^\sharr $ is automatically a $\cup$-homomorphism, hence a 
$ \subseteq $-homomorphism.
If
$\varphi$  is injective, then 
$\varphi ^\sharr $ is an embedding with respect to
both $\cup$  and $ \subseteq $. 

If  $\varphi$  is continuous and $ a \sqsubseteq b$, that is, 
$a \subseteq  \K b$, we have 
$\varphi ^\sharr  (a) \subseteq \varphi ^\sharr (  \K b) 
 \subseteq   \K  \varphi^ \sharr (b)$, since $\varphi ^\sharr$
is a $ \subseteq $-homomorphism.    
Hence $ \varphi ^\sharr  ( a) \sqsubseteq \varphi ^\sharr  (b)$,
thus $\varphi ^\sharr  $ is a $ \sqsubseteq $-homomorphism.
We have proved (1) $\Rightarrow $  (2). 
The implication (2) $\Rightarrow $  (3) is immediate from the 
fact that a semilattice homomorphism
is an order preserving map for the induced order \cite[p.\ 30]{G}. 

To prove (3) $\Rightarrow $  (1),  suppose that  
$\varphi ^\sharr  $ is a $ \sqsubseteq $-homomorphism, that is,
$a \sqsubseteq b$ implies 
$ \varphi ^\sharr  ( a) \sqsubseteq \varphi ^\sharr  (b)$,
for every $a,b \subseteq X$. In particular, we can take
$a= \K b$ and, since   
$ \K b \subseteq  \K b$, we get
$  \K b \sqsubseteq b$, hence 
$ \varphi ^\sharr  (  \K b) \sqsubseteq \varphi ^\sharr  (b)$, thus
$ \varphi ^\sharr  (  \K b) \subseteq  \K \varphi ^\sharr  (b)$.
Hence $\varphi$  is continuous. 

To prove the last statement, first observe that
the following is a chain of equivalent conditions, for $Z$
a 
closure space over some subset of $Y$.
  \begin{enumerate}[(i)]    
\item   
 $Z $ is a subspace of $ Y$, 
\item
the corresponding closure operations satisfy
$ \K _Z d = Z \cap  \K _Y d$, for all  
 $d \subseteq Z$.
\item
for all  $c,d \subseteq Z$,
  $c \subseteq  \K _Z d$
if and only if $c \subseteq  \K _Y d$.
\item
for all  $c,d \subseteq Z$,  $c \sqsubseteq  _Z d$
if and only if $c \sqsubseteq _Y d$.
 \end{enumerate}

  Thus if $\varphi: X \to Y$
is an embedding
of 
closure spaces,
$Z= \varphi ^\sharr (X)$ 
and $a,b \subseteq X$, then
$\varphi ^\sharr (a) \sqsubseteq _Y \varphi ^\sharr (b)$
if and only if 
$\varphi ^\sharr (a) \sqsubseteq _Z \varphi ^\sharr (b)$,
if and only if 
$ a \sqsubseteq_X b $,
since $\varphi$  induces a homeomorphism
from $X$ onto  $Z= \varphi ^\sharr (X)$.

 Conversely, if $\varphi ^\sharr : \mathcal P(X) \to \mathcal P(Y) $
is an embedding of specialization semilattices
from   $\mathbf S(X)$ to $\mathbf S(Y)$,
then $\varphi : X \to Y $ is injective.
Endow  $Z= \varphi ^\sharr (X)$ with the closure induced by 
the closure on $X$ through $\varphi$,
thus   $\varphi ^\sharr (a) \sqsubseteq _Z \varphi ^\sharr (b)$
if and only if 
$ a \sqsubseteq_X b $.
This is also equivalent to
 $\varphi ^\sharr (a) \sqsubseteq _Y \varphi ^\sharr (b)$,
since $\varphi ^\sharr  $
is an embedding of specialization semilattices.
Since $\varphi ^\sharr$ is surjective from 
$X$ to $Z$, then, for every $c,d \subseteq  Z$, there are 
$a,b \subseteq  X$ such that  
$c= \varphi ^\sharr (a) $ and $ d=  \varphi ^\sharr (b)$.
By the equivalence of (iv) and (i), $Z$ is a subspace of $Y$
and this means precisely that  $\varphi$  is an embedding
of closure spaces.
\end{proof} 

\begin{remark} \labbel{otoh}    
On the other hand, as we mentioned before, 
if $\varphi$  is continuous, it is 
not necessarily the case that 
$\varphi ^\sharr ( \K b) \supseteq   \K  \varphi^ \sharr (b) $.
In fact, equality holds if and only if 
$\varphi$  is a closed map \cite[Exercise 1.4.C]{En}. 
If $\psi$ is a function 
between two sets with some unary operation $ \K $,
we say that
$\psi$ is a \emph{homomorphism with respect to $ \K $},
a \emph{$\K$-homomorphism}, for short, 
if $\psi ( \K b) =  \K  \psi (b)$,     
for all elements $b$ in the domain of $\psi$.
In this terminology,
$\varphi$  being a continuous function between
two topological spaces does not entail 
$\varphi^\sharr$ being a   homomorphism with respect to $ \K $.
\arxiv{

Notice that
 it is not the case
that to every function $\psi: \mathcal P(X) \to \mathcal P(Y) $
there is an associated function
$\varphi : X \to Y$  such  that 
$\psi= \varphi ^\sharr$,    
let alone continuity and the notion of homomorphism. } 
 \end{remark}

It follows from Proposition \ref{propcont} 
that a function $\varphi$    between two topological 
\arxiv{or closure }
spaces is
continuous if and only if 
$\varphi^\sharr$ preserves the equivalence relation 
$\equiv$  among  
subsets defined by $x \equiv y$ if $x$ and $y$
have the same closure (though $\varphi^\sharr$  does not
necessarily preserve the closure itself!). 
\arxiv{The author believes that it is not yet completely clear why 
 the most used notion of
morphism in general topology is a continuous function. As stressed 
by Kuratowski, continuous functions are exactly functions which
preserve the adherence relation (and 
the analogue for specialization here is just an immediate
generalization of this observation). However, if we argue this way, many
scholars suggest that proximity is a much more natural notion than topology
\cite{DC}. 
Could the above observation (that continuous functions are
exactly those functions preserving the condition of 
 having equal closure) contribute to explain
 the usefulness of continuous functions? }

Many of the above definitions, ideas and arguments
apply to a general setting in which the underlying structure
is just a poset.

\begin{definition} \labbel{closop}
If  $(P, \leq)$ is a poset,
a \emph{closure operation} 
is an isotone, inflationary and idempotent unary operation
$ \K $  on $P$. 
In the above situation the triple
 $(P, \leq,  \K )$ shall be called a \emph{closure poset}. 
If $\leq$ is associated to some semilattice operation 
$\vee$, we shall say that
  $(P, \vee,  \K )$ is a \emph{closure semilattice}.
Again, see \cite[Section 3.1]{E}
for further details and an equivalent characterization
in terms of ``closed'' elements. 
\arxiv{See also Remark \ref{sc}(b) below. 

}
 If $(P, \leq)$ is a poset
and $ \K $ is a closure operation
on $P$, the 
 \emph{associated specialization poset}
is the model $( P, \leq , {\sqsubseteq} )$,
where $ \sqsubseteq $ is  as in Definition \ref{asso},
namely, $ a \sqsubseteq b$
if  $a \leq  \K b$,
for $a, b \in P$.
If in addition  
$P$ is a join-semilattice, we will
speak of the \emph{associated 
specialization semilattice} $\mathbf S(P)$. 

A \emph{homomorphism} between two 
closure posets (semilattices)
is an order-preserving map (a semilattice homomorphism) which is also a 
homomorphism with respect to 
$ \K $, as defined in Remark \ref{otoh}.
Any homomorphism between two closure posets 
is a homomorphism between the 
associated 
specialization semilattices; however, the converse
does not necessarily hold, as already seen 
in the special case  of topological spaces. See
Remark \ref{otoh} 

If we introduce a notion of 
``continuity'' between closure posets,
a generalization of Proposition \ref{propcont}
holds with the same proof.  If $\mathbf P$ and  $\mathbf Q$
are closure posets (semilattices), we say 
that a function $\psi: P \to Q$
is \emph{continuous} if $\psi$ is  order-preserving
(a semilattice homomorphism) and  moreover 
  $\psi ( \K _P  a) \leq  \K _Q \psi (a)$,
for every $ a \in P$.  Continuity between closure posets
has been studied in \cite{Sla}, where the following proposition
has been essentially proved. 
\end{definition}

\begin{corollary} \labbel{corr} 
An order preserving function  (semilattice
homomorphism) $\psi$ between two closure posets
(semilattices) is continuous if and only if 
$\psi$ is a  homomorphism between the 
associated 
specialization posets (semilattices).
 \end{corollary}

\arxiv{ 
It is probably interesting to observe that no assumption on $ \K $ 
is necessary in order to prove Corollary \ref{corr}.
In more detail,   
if $\mathbf A$ and  $\mathbf B$
are ordered models for a language $\mathcal L$ 
with a unary function 
$ \K $, let us say 
that a function $\psi: A \to B$
is \emph{continuous} if $\psi$ is a homomorphism
for the $\mathcal L \setminus \K $ reduct,
  (in particular, $\psi$ is order-preserving)
and  moreover 
  $\psi ( \K _A  a) \leq  \K _B \psi (a)$,
for every $ a \in A$.   
Define $ \sqsubseteq $  by 
$ a\sqsubseteq b $ if 
$ a \leq \K  b $.

\begin{corollary} \labbel{corrr} 
A function $\psi$ between two models
$\mathbf A$ and $\mathbf B$ as above
 is continuous if and only if 
$\psi$ is a  homomorphism 
in the language $ (\mathcal L \setminus \K ) \cup \sqsubseteq  $,
where $ \sqsubseteq $ is defined as above. 
 \end{corollary}

Actually, we do not need the assumption that $\leq$   is an
order relation, to give the above definition  of continuity.
Corollary \ref{corrr} then holds just under the assumption that $\leq$  is a 
reflexive relation.
  }

\section{Specialization semilattices and posets} \labbel{specsec} 

\arxiv{Recall that a \emph{preorder} 
on some set $P$ is a binary
reflexive and transitive relation on $P$. 
Some authors use the expression
\emph{quasiorder} in place of preorder.
A \emph{partial order}, or simply \emph{order}, or \emph{ordering}
or \emph{order relation} is an antisymmetric preorder.
A \emph{partially ordered set},
for short, \emph{poset}
is a set   endowed with a partial order.
 
Recall that (algebraically) a  \emph{semilattice} is a set $S$ 
together with a binary operation
$\vee$ which is commutative, associative and idempotent.
Recall that semilattices  here will be always  intended as join-semilattices,
in the sense that a semilattice $(S, \vee)$
induces the partial order $\leq$ on $S$
defined by 
$a \leq b$ if and only if $a \vee b =b$.
Meet semilattices---not considered here---algebraically
are defined in the same way, but conventionally
are assumed to induce the reverse order.
When we mention a partial order in reference to some 
semilattice, we will always mean the order $\leq$ 
introduced above. 
 }

\begin{definition} \labbel{spsem}
(a)
A \emph{specialization poset}
$\mathbf S$ is a structure
$(S, \leq, {\sqsubseteq} ) $ 
such that 
$(S, \leq)$
is a poset and
$ \sqsubseteq $ is a binary relation on $S$ satisfying 
\begin{align}
\labbel{s1}    \tag{S1}
 & a \leq b  \Rightarrow  a \sqsubseteq b, 
&& \text{($\leq$ is finer than $ \sqsubseteq $)} 
\\
\labbel{s2}    \tag{S2}
& a \sqsubseteq b \ \&\   b \sqsubseteq c \Rightarrow 
 a \sqsubseteq c, 
&& \text{($ \sqsubseteq $ is transitive)} 
 \end{align}    
for all elements $a,b, c \in S$.

(b)  
A \emph{specialization semilattice}
$\mathbf S$ is a triple 
$(S, \vee, {\sqsubseteq} ) $
such that 
$(S, \vee)$
is a semilattice, $\mathbf S$ satisfies 
\eqref{s1} - \eqref{s2},  where 
$\leq$ is the order induced by $\vee$, and
\begin{align}
\labbel{s3}    \tag{S3}
 &a \sqsubseteq b  \ \&\  a_1 \sqsubseteq b
 \Rightarrow 
 a \vee a_1 \sqsubseteq b, 
&& \text{($ \sqsubseteq $ respects joins on the $1$st comp.)} 
\end{align}
 for all elements $a,b,   a_1  \in S$.

In both cases, the relation $ \sqsubseteq $ shall be called
a \emph{specialization}.  
\end{definition}   

In words, a specialization poset  
is a  poset  endowed   with an additional preorder 
$ \sqsubseteq $ which is coarser (that is, larger)
than the poset order. See Remark \ref{conseq}(c) below.
In the case of specialization semilattices we also ask for the compatibility
condition \eqref{s3}.
An infinitary version of \eqref{s3}
for complete lattices (resp., complete join-semilattices)
with a coarser preorder has been considered
in \cite[Subsection 3.1]{GT},
\cite[p.\ 232]{Now}, \cite[Definition 2.1]{Sla}.
We present below the very elementary example of 
a join semilattice with  the structure of 
a specialization poset, but which is not a specialization semilattice.

 \begin{example} \labbel{diamond}
Let $S= \{ 0, a, b, 1\} $,
with the partial order $\leq$  
given by 
$0 < a < 1$
and  $0 < b < 1$.
Let $1 \sqsubseteq 1$,
 $ x \sqsubseteq y$ and
$  x \sqsubseteq 1$,
for all $x,y \in \{ 0,a,b \} $ 
and let no other 
$ \sqsubseteq $-relation hold.

  \begin{center}
\begin{tikzpicture}[-,>=stealth',auto,node distance=3cm,
thick,main node/.style={circle,scale=0.3,draw},
mnode/.style={scale=0.3}]
\node[main node,label={south:$0$}] (0){};
\node[main node,label={west:$a$}] (a) [above left of=0]{};
\node[main node,label={east:$b$}] (b) [above right of=0]{};
\node[main node,label={north:$1$}] (1) [above right of=a]{};
\draw (0) -- (a)-- (1);
\draw (0) -- (b)-- (1);
\end{tikzpicture}
\end{center}

Then $(S, \leq , {\sqsubseteq} )$
is a specialization poset.
Moreover, $\leq$ induces  
a semilattice operation $\vee$ on $S$, 
but \eqref{s3} fails in  
$(S, \vee, {\sqsubseteq} )$,
since $ a \sqsubseteq 0$,
$  b \sqsubseteq 0$ but
$a \vee b = 1 \centernot \sqsubseteq 0$. 
Hence $(S, \vee, {\sqsubseteq} )$ 
is not a specialization semilattice.
 \end{example}

\begin{remark} \labbel{si!}   
The structures introduced in 
Definition \ref{asso} and, more generally, in Definition \ref{closop} 
are easily seen to be specialization semilattices and posets
according to the preceding definition.  
\end{remark}
 
\begin{definition} \labbel{osr}    
If $\mathbf S = (S, \vee, {\sqsubseteq} ) $
is a specialization semilattice, then
$\mathbf R (\mathbf S) = (S, \leq, {\sqsubseteq} ) $
is a specialization poset, which shall be called the
\emph{order-specializa\-tion-reduct} of  $\mathbf S$. 
Strictly speaking, this  is not a reduct 
of $\mathbf S$ 
in the formal sense 
(it is a reduct of
$ (S, \vee, \leq, {\sqsubseteq} ) $) but 
we hope the terminology is sufficiently clear and intuitive. 

In particular, everything we will say about
specialization posets will apply
to specialization semilattices, too. 
\end{definition}

\begin{remarks} \labbel{conseq}
(a)
From \eqref{s1} one immediately gets 
\begin{align} 
\labbel{s4}    \tag{S4} 
 & a  \sqsubseteq a, 
\quad\quad\quad  \text{($ \sqsubseteq $ is reflexive)} 
\\
\intertext{while from \eqref{s1} - \eqref{s2} one gets}  
\labbel{s5}    \tag{S5}
 a \sqsubseteq b \ \&\  b \leq c
& \Rightarrow 
 a \sqsubseteq c, 
\quad \text{($ \sqsubseteq $ is weakly transitive on the right)} 
\\
\labbel{s6}    \tag{S6}
 a \leq b \ \&\ 
b \sqsubseteq c & \Rightarrow 
 a \sqsubseteq c,
\quad \text{($ \sqsubseteq $ is weakly transitive on the left)}
\\
\intertext{for all elements $a,b, c \in S$.
Thus every specialization poset
satisfies \eqref{s4} - \eqref{s6}. 
From $b,b_1 \leq b \vee b_1$, \eqref{s5}
and \eqref{s3}  it follows}
\labbel{s7}    \tag{S7}
 a \sqsubseteq b  \ \&\  a_1 \sqsubseteq b_1
 & \Rightarrow 
 a \vee a_1 \sqsubseteq b \vee b_1, 
\quad \text{($\vee$ preserves $ \sqsubseteq $)}
\\
\intertext{in particular, by \eqref{s4} and taking $a_1=b$
in \eqref{s3}, respectively, $a_1=b_1$
in \eqref{s7}, we get}
\labbel{s8}    \tag{S8}
 a \sqsubseteq b  
 & \Rightarrow 
 a \vee b \sqsubseteq b , 
\\
\labbel{s9}    \tag{S9}
 a \sqsubseteq b  
 & \Rightarrow 
 a \vee a_1 \sqsubseteq b \vee a_1, 
\end{align}
thus \eqref{s4} - \eqref{s9} hold in every
specialization semilattice. 

(b) From the above remarks  we get alternative axiomatizations.
For example, a poset with a further relation $ \sqsubseteq $
is a specialization poset if and only if 
\eqref{s2}, \eqref{s4} and \eqref{s5} hold
(equivalently, \eqref{s5} can be replaced by \eqref{s6}).

 (c) It follows from 
\eqref{s4} and  \eqref{s2} that
if $\mathbf S$ is a specialization poset 
(semilattice), then  $ \sqsubseteq $
is a preorder. However,
$ \sqsubseteq $
is not required to be antisymmetric,
hence it is not necessarily an order.
In particular, whenever we speak of \emph{meets} and
\emph{joins}, these notions shall be always
interpreted as relative to $\leq$.  

(d)
Notice that properties 
\eqref{s1} - \eqref{s9} are all expressible as    
first-order universal sentences: just prefix each formula
by an appropriate string of universal quantifiers.
Actually 
\eqref{s1} - \eqref{s9} are Horn formulas \cite[Section 9.1]{H} 
and this implies that they are preserved under direct
products, a fact which can also be easily verified directly.
\end{remarks}   

\arxiv{
\subsection{Homomorphisms and embeddings} \labbel{homomorphisms} 
For the sake of completeness, we give the explicit
definitions of homomorphisms and embeddings
between specialization posets and  semilattices.
The notions correspond exactly to the standard model-theoretical
notions \cite[Section 1.2]{H}. 
Some authors use the expression
\emph{isomorphic embedding} \cite[Section 1.3]{CK} 
for what we call  simply  an embedding.

If $(P, \leq_P)$ and $(Q, \leq_Q)$
are posets and $\varphi: P \to Q$,
then $\varphi$  is an \emph{order-preserving map},
or an \emph{ordermorphism}, or simply
a \emph{homomorphism} when the context is clear,        
if $a \leq_P b$ implies $\varphi(a) \leq _Q \varphi (b)$,
for every $a,b \in P$.
An order-preserving map is an \emph{order-embedding},
or simply an \emph{embedding} if in addition
 $\varphi(a) \leq _Q \varphi (b)$ implies $a \leq_P b$,
for every $a,b \in P$. Notice that an order-embedding
is necessarily injective.

If $(S, \vee_S )$ and $(T, \vee_T)$
are semilattices, a \emph{(semilattice) homomorphism} 
is a function $\varphi:S \to T$ such that 
$\varphi(a \vee _S b) = \varphi (a) \vee_T \varphi(b) $,
for every $a,b \in S$. A \emph{(semilattice) embedding}
is an injective homomorphism.
Notice  that a semilattice  homomorphism (embedding) 
between two semilattices induces an order-preserving map
(an embedding) for the corresponding 
``order-reducts''.

If $(S, \leq_S, \sqsubseteq_S )$ and $(T, \leq_T, \sqsubseteq _T)$
are specialization posets, $\varphi$  is
a \emph{homomorphism (of specialization posets)}
if $\varphi$ is an order-preserving map from
$(S, \leq_S )$ to $(T, \leq_T)$ and moreover
\begin{equation}\labbel{M} \tag{M}      
  \text{$a \sqsubseteq _S b$ implies $\varphi(a) \sqsubseteq  _T \varphi (b)$,
for every $a,b \in S$.}
 \end{equation}
A homomorphism of specialization posets
is an \emph{embedding}
if it  is an order-embedding and moreover
 \begin{equation}\labbel{E} \tag{E} 
\text{    $\varphi(a) \sqsubseteq  _T \varphi (b)$ implies
$a \sqsubseteq _S b$,
for every $a,b \in S$.} 
\end{equation} 

If $(S, \vee_S, \sqsubseteq_S )$ and $(T, \vee_T, \sqsubseteq _T)$
are specialization semilattices, a 
 \emph{homomorphism (of specialization semilattices)}
is a semilattice homomorphism  satisfying \eqref{M}. 
A homomorphism of specialization semilattices
is an \emph{embedding} if it is injective and satisfies \eqref{E}. 

As already mentioned, for models with a unary operation $ \K $, 
e.g.\   closure posets as introduced in Definition \ref{closop},
a homomorphism is supposed to satisfy also 
$\varphi( \K _S a) =  \K _T \varphi (a)$.

In the above definitions we have distinguished, say,
the operation $ \vee_S$ on $S$ from the operation 
$ \vee_T$ on $T$ by adding the corresponding subscripts.
When no risk of ambiguity might occur, we will drop the subscripts.
To be more accurate, we should have made the distinction
 between  symbols
and their interpretations \cite{H,CK}.
Here we do not need to make the distinction explicit,
hence we will proceed quite informally.

Notice that in algebra (namely, when dealing only with
operations) there is no distinction between
 injective homomorphism  and embeddings.
On the other hand, in topology (and in model theory
when predicates, i.e. relations, are present)
injective  continuous functions (or homomorphisms) 
are not necessarily embeddings.
For example, if $X$, $Y$ are topological spaces, $\varphi: X \to Y$
and
 $\varphi$  is continuous and injective, it is not
necessarily the case that $X$ is homeomorphic to a subspace of
$Y$. 
 }

\section{Further examples} \labbel{exasub} 

In this section we describe more examples
of specialization semilattices and posets. 
The material in this section is
presented only as a further motivation;
the section is not necessary for  the remaining
parts of the paper.

\smallskip 

\emph{1. The specialization preorder.} 
Given a topological space $X$, the relation 
$ x \sqsubseteq y$ between \emph{points} of $X$,
defined by $x \in  \K  \{ y \} $ is called the 
\emph{specialization preorder}. 
It has
 various applications, among others, to algebraic geometry
\cite[Ex. 3.17e]{Ha} and to domain theory \cite{GHKLMS}.
To make the example fit with our abstract definition 
of a specialization poset, we let the underlying order on $X$ 
be the trivial order in which no two distinct  points are comparable. 

In this note we  have simply considered the extension of the
specialization  preorder to all \emph{subsets} of $X$. 
 
\arxiv{For algebraic closure spaces, 
a multiary generalization of the specialization preorder
 has been considered 
by A. Pasini \cite{Pa}, who also
noticed  that in this situation  homomorphisms correspond
exactly to continuity. However, a topology
(which usually is not algebraic)
cannot be generally retrieved by the closure of finite subsets. }

\emph{2. Structures with two comparable binary relations.}
Various structures with a pair of comparable preorders have been 
considered in the literature,
for example, in connection with semantics 
of  modal logics \cite{BHl,BH,FM,GKWZ},
domain theory \cite[Sect.\ I-1]{GHKLMS},  
tolerance spaces \cite{PW},  measures \cite{L},
 representations of lattices \cite{Hol,San,U}, 
and even 
abstract foundations of 
general relativity \cite{KP,Pp}. 
Strictly speaking, in some cases the structures we have referred to are not
examples of specialization posets. For example, formally,
 in
Fairtlough and Mendler frames \cite{BH,FM}
both relations are only preorders.
 However,  in many situations,
one can do with a  partial order
\cite[Proposition 4.5]{BHl}, \cite{BH}.
 
A  \emph{tolerance space}
is a set $X$ together with a symmetric and reflexive
relation $\tau$ on $X$.
\arxiv{According to \cite{PN}, tolerance spaces have been introduced by 
E. C. Zeeman in the 1960s, but
A. B. Sossinsky observed that the main
idea underlying tolerance spaces comes from
J. H. Poincar\'e. }
If $X$ is a tolerance space, then,
 for $a,b \subseteq X$, let $ a \sqsubseteq_ \tau  b $ if
$\tau(a) \subseteq \tau(b)$, where
 $\tau(a)=
\{ x\in X \mid y \mathrel { \tau  } x,
\text{ for some } y \in a \} $.
In the present 
terminology, $\cup$ and $ \sqsubseteq _ \tau  $ define
the structure of a specialization semilattice 
on  $\mathcal P(X)$.
The above notions have found applications to image
analysis  and other information systems
\cite{PW}. 

If $\mu $ is a measure defined on some set $S$ of subsets
of $X$,  let  
$a \sqsubseteq_ \mu b$
if $\mu (a) \leq \mu (b)$, for $a,b \in S$.
Then subset inclusion and $ \sqsubseteq _\mu$ 
provide $S$ with the structure of a specialization poset and,
if $\mu $ is 
two-valued, of a specialization semilattice. 
   The relation $ \sqsubseteq_ \mu  $ has been widely studied
in connection with foundations of probability and with purported
economical applications \cite{L}.  
\arxiv{See Example \ref{meas}
for more details. }

In connection with \cite{San}
we mention that if $\mathbf P$ is a poset and 
$\mathcal P_f (P)$ is the set of all finite subsets of 
$P$, then    $(\mathcal P_f (P), {\cup}, {\ll})$
is a specialization semilattice, where $\ll$
is the \emph{refinement}  relation on $\mathcal P_f (P)$ defined by
$X \ll Y$ if, for every $x \in X$,
there is $y \in Y$ such that 
$x \leq_{\mathbf P} y$. 
See also \cite[p. 283]{G}.

\arxiv{\emph{2b.Causal spaces}
have been
introduced by E. H. Kronheimer and R. Penrose in \cite{KP} 
 in connection with 
abstract foundations of 
general relativity. 
Causal spaces are sets with two comparable partial orders
satisfying a further  coherence condition
connecting the two relations\footnote{In \cite{KP}
 a further relation $\sharr$  has been considered; 
however, $\sharr$ is definable by means of the two posets.
Formally, in \cite{KP} one poset is considered as an antireflexive relation,
but the theory  is biinterpretable with a theory 
with two (reflexive) posets, one finer than the other.};
in particular, causal spaces are specialization posets.
See Remark \ref{causs} below for further comments.} 

\emph{3. Specializations induced by a quotient.}
If $\mathbf S$, $\mathbf T$  are semilattices
\arxiv{(posets)}, 
$\varphi: \mathbf S \to \mathbf T $  is a homomorphism
\arxiv{(an order preserving map)}, 
 and we let $ a \sqsubseteq_ \varphi  b$ in $S$ 
if $\varphi(a) \leq \varphi (b)$
in $\mathbf T$, then 
$\mathbf S$ 
is endowed with the structure of a specialization 
semilattice \arxiv{(poset)}. 

In particular, if $\mathbf B$
is a Boolean algebra and $\mathcal I$
is an ideal on $\mathbf B$, 
then $(B,\vee, {\sqsubseteq} )$
is a specialization semilattice, 
where $a \sqsubseteq b$
if $a/\mathcal I \leq b / \mathcal I$
in  $\mathbf B/\mathcal I$.
Many examples of such structures have been widely studied;
we mention just one.
 If $B= \mathcal P ( \mathbb N )$ and
$\mathcal I$ is the set of all finite subsets of $ \mathbb N$, then
$ \sqsubseteq $ is inclusion mod finite.
The notion  has applications to descriptive \cite{MN} and combinatorial set theory
\cite{Bl,Fa}, topology \cite{DH}, model theory \cite{MN},
among others. Note that here there is no underlying notion of
``closure'': for every subset $y$ of $\mathbb N$, there are many
larger subsets $x$ such that $ x \sqsubseteq y$, that is, $x \setminus y$
is finite. However, there is no largest such $x$.

\arxiv{Note that the specialization semilattice (poset)
associated to a tolerance (to a measure), as described
in the previous subsection, are particular examples 
of the present quotient construction. For example, in the case
of tolerances, the function $\varphi$   which assigns $\tau(a)$  
to $a$  is a semilattice homomorphism  from 
$(\mathcal P(X), \cup)$  to $(\mathcal P(X), \cup)$,
thus $ \sqsubseteq_ \varphi $ here is the same as $ \sqsubseteq_ \tau $  in 
2.} 
 
The definition of $a \sqsubseteq_ \varphi b$ above,
for $\varphi$  an arbitrary semilattice homomorphism,
is very general; we will show in a sequel  to the present
note that every specialization semilattice  
can be constructed this way. 
\arxiv{ See \cite{mttcong} or \cite[Theorem 5.6]{mttmult}. 

Moreover, we get an essentially  equivalent 
definition of a specialization semilattice 
 if we  consider a semilattice together with a congruence. 
See Definition \ref{cong}.
For short, a specialization semilattice can be seen 
as a substructure of (some structure associated to) a topological spaces,
but can equivalently be seen also as a semilattice together 
with an onto homomorphism, equivalently, a 
semilattice together with a congruence. This correspondence
can be also interpreted in a categorical setting. We refer again to
\cite{mttcong} for more details.}

\emph{4. Abstract consequence relations on posets.} 
Abstract consequence relations  have been first introduced
by A. Tarski (see \cite{T}) in a fashion
slightly different from the modern treatment \cite[Note on p. 23]{W},
\cite[Section 1]{Jan}.
In the restricted classical sense, an \emph{abstract consequence relation}
 is a binary relation
$\vdash$ between sets of formulas and formulas of 
a formal language. The intended meaning of 
$\Gamma \vdash \sigma $ is that 
$\sigma$ is deducible from $\Gamma$
in some deduction system fixed in advance.
Consequence relations can be introduced
abstractly and provide an equivalent formulation
for the notion of a  \emph{consequence operation}
or a \emph{closure operation} \cite{GT,Jan}.
To some consequence relation $\vdash$
one associates the   closure operation $C$ 
sending a set $\Gamma$ of formulas to the set 
$C(\Gamma)$ of all the formulas deducible from $\Gamma$.
The correspondence with topology is patent:
a consequence operation is interdefinable with a
  consequence relation 
exactly in the same way as 
topological closure is interdefinable with 
the adherence relation.

If one abstracts from single formulas ``being singletons'',
one can  introduce a binary relation
$ \Gamma \vdashcup  \Sigma $  between
sets of formulas, whose intended meaning
 is 
``everything in $\Sigma$ is deducible from $\Gamma$.''\footnote{Most authors
use $ \Gamma  \vdash \Sigma $ to mean that 
``at least one formula  in $\Sigma$ is deducible from $\Gamma$'':
this is the reason for our notation.
See e.g.\  \cite{BBI,Ie,Je}; see \cite{fiore,No} 
  for a 
comparison between the two situations.}
 Notice the resemblance with  \eqref{1}
(of course,  considering the converse relation).
This  is not 
love of abstraction for its own sake;
the idea provides the possibility of introducing 
 consequence relations in the setting of
arbitrary complete lattices, or even
posets, with deep and important applications 
to algebraic logic, in particular, concerning equivalence and algebraizability
 \cite{GT}.

 Note that a remark 
similar to Remark \ref{anchclos} applies to the present situation.
If some deduction system with consequence operation
$C_2$ extends a system with consequence $C_1$,
then, in general, $C_1( \Gamma ) $ is smaller than
 $C_2(\Gamma ) \cap \Sigma  $, where $\Sigma$ is the set
 of  sentences of the first system \cite[1.7.6]{W}.
Equality holds only in the case of conservative extensions.
On the other hand, the ``specialization''  relation 
$\vdashcup$ is preserved in passing from a 
deduction system  to an extension. 

\arxiv{It is claimed in \cite[Section 3]{GT} that a large part of the results
and definitions in that section apply to arbitrary partially ordered sets
in place of complete lattices.
In the setting of partially ordered sets
 Conditions (1) - (2)
in \cite[Subsection 3.1]{GT}  
correspond to the definition of a 
specialization poset  as presented here\footnote{again,
  considering the converse of the relation $\vdashcup$.}.
In the case of join semilattices we get specialization semilattices if we add 
the finitary version of
Condition 3  from \cite[Subsection 3.1]{GT}. If we 
add the full Condition 3, we get principal specialization semilattices; see
Definition \ref{supcomp} below. 
Compare also \cite[Definition 2.1]{Sla}. }

\emph{5. Recursive sets of formulas.}
As just mentioned, given a  
deduction system, it is natural to define
a relation $ \Gamma  \vdashcup \Sigma $
between sets of formulas, meaning 
that everything in $\Sigma$ is deducible from $\Gamma$.
In the framework of arbitrary sets of formulas, 
the definition of  $   \vdashcup $ 
is interchangeable with the definition of 
the consequence operation $C$
which assigns to a set $\Gamma$ the set
$C(\Gamma)$ of all the formulas deducible from $\Gamma$.
In fact, from the closure operation $C$ we can recover 
$\vdashcup$ by setting  
$ \Gamma  \vdashcup \Sigma $ if 
$\Sigma \subseteq C(\Gamma )$.
Notice again the similarity with \eqref{1}. 
Thus, when dealing with  arbitrary sets of formulas,
it is generally irrelevant whether we deal with 
$\vdashcup$ or $C$.  

Suppose now that we want to deal only with finite sets of formulas, 
or, possibly, with recursive sets of formulas.
The assumption makes good sense,
since only recursive sets of formulas can be effectively described.
In this framework the  approach 
using  consequence operations is not equivalent to the approach
using consequence relations; actually, the former 
is not viable, since
in general the set
$C(\Gamma)$ of the formulas deducible from $\Gamma$
is not recursive, even when $\Gamma$ is.
This is an example showing that ``specialization''
is more apt than ``closure'' in certain situations. 

\arxiv{\emph{5b. Finite pieces of information.} 
The above example can be immediately reformulated in a 
much more practical setting.
 Suppose that we are dealing with finite pieces of information,
such as, e.g., information that can be stored in some database.
Suppose that we also have some constraints on our data, so that
from some information we may obtain more information.
Given a set $\Gamma$ of such information, the set 
$C(\Gamma)$ of all the informations which can be obtained from $\Gamma$  
might be too large to be stored, hence 
$C(\Gamma)$ 
 is generally
 an object inappropriate in our framework.
The relation
$ \Gamma  \vdashcup \Sigma $,
meaning that all the informations in $\Sigma $ 
can be obtained from the informations in $\Gamma$ 
is surely more concrete and manageable.

If we consider also the relation of plain containment (union),
then the set of the pieces of information under consideration
becomes a specialization poset (semilattice) according to our definitions.

More formal details are presented in the next example.
 }

\emph{6. Complete implicational systems.}
The axioms for a specialization semilattice 
are equivalent to the properties characterizing a
\emph{complete implicational system} 
on the power set of a finite set 
\cite[Section 7.4, in particular Theorem 7.65]{CLM},
\cite{CM}, considering the reverse relation.
Complete implicational systems, variously called \emph{complete families of
functional dependencies}, \emph{entail relations},
\emph{closed families of implications} or \emph{implicational 
theories} \cite{CM} have found applications 
in  relational data bases,
data analysis,   
 artificial intelligence and  mathematics
of social sciences. See \cite{CLM,CM} for references.
The notion of specialization presented here is much more general,
since it applies to any, possibly infinite, semilattice.
For example, strictly formally, the relation
$ \Gamma  \vdashcup \Sigma $ introduced in 5.\ above
is not the implication of a 
complete implicational system, since it is defined on an infinite set.

\arxiv{In particular, specialization semilattices
seem to be a promising generalization of 
complete implicational systems. The axioms for CIS
are the same, but, as far as we
know, CIS have been considered only on the power set
of  a finite set. As just mentioned, CIS have many applications.
 Infinitary structures similar to CIS 
have been considered (for different purposes), but only---still, as
far as we know---in the case of complete semilattices. Dealing with (finitary)
semilattices seems to have some advantages. 
The above line of inquiry is only hinted here,
but we believe that it is potentially interesting.

In more detail, an infinitary generalization of complete implicational
systems (CIS)  is connected to 
contexts in which potentially infinite objects can be considered.
Intuitively, an instance of an implication in a CIS is 
an assertion of the form ``everything having such and 
such characteristics has also such and such other characteristics''
\cite[p.\ X]{CLM}. 
On the set of all finite subsets of natural numbers, we get a 
specialization semilattice if we let $x \sqsubseteq y$  mean
``every natural number which can be written
as a finite product of numbers in $x$, possibly with repetitions,
can be written
as a finite product of numbers in $y$, possibly with repetitions''.
Such a relation cannot be described on a (finite) CIS.
Underlying the above relation, of course, 
there is a closure operation $ \K $ defined by
 taking $ \K  y $ to be the set of all numbers
which can be written as a product of numbers in $y$.
However, $ \K  y $ is generally infinite, hence, 
if we want to remain in the realm of finite objects,
an implicational or ``specialization''  relation of the form $x \sqsubseteq  y  $ is
admissible, while a set as $ \K  y $  is not.
Note that  this example is very similar in spirit to the example in 
5.\ above. 

\smallskip 

As we have showed, there are many situations in which
specialization appears to be a more natural notion than closure.
Moreover, as already remarked, specialization frequently
behaves better with respect to  homomorphisms. This has been explained at
length in Sections \ref{intro}  and \ref{motiv} for the case of topological spaces. Other examples 
have been presented in Remark \ref{anchclos}  and in 4. above:
 see the remark about conservative extensions.
In general, the main point  is that
the closure operation is not preserved on submodels 
of, say, closure posets. 
On the other hand, if $\mathbf S$ is a specialization poset (semilattice),
then the specialization is preserved on every subset (sublattice) of $\mathbf S$.
In more detail, if $X$ is a closure space,
 then the closure operation naturally induced on 
some subset $Y$ of $X$ acts differently on the subsets of $Y$.
In other words, we cannot maintain the same closure operation on $Y$,
unless $Y$ is closed in $X$. The same applies to topological spaces.

As a way of example, if $Y {\subseteq} X$ are topological spaces,
then $( \mathcal P(Y), {\subseteq} , \K _Y )$
is not necessarily a submodel of $( \mathcal P(X), {\subseteq} , \K  )$,
while $( \mathcal P(Y), {\subseteq} , \sqsubseteq  )$
is indeed a submodel of $( \mathcal P(X), {\subseteq} , \sqsubseteq   )$.
  }

\smallskip 

In the present section we have limited ourselves to examples
which match exactly the definition of a specialization semilattice 
or poset. Needless to say,  many variations or extensions
of the above notions have been considered in the literature.
Relational and algebraic structures associated or related to topologies
\cite{Ed,GHKLMS,Iv,J,MP,PP}
 and to
modal logics \cite{CPZ,FM,Fr,GM,Jan,Je,W}
have a long history;  surveys with many examples 
and discussing the interconnections can be found in 
\cite{vBB,BH,Es,GKWZ,GH,Gr}.

Papers very similar in spirit with the present one include 
\cite{DIM,Iv}.
In particular, \cite{Iv} presents 
further motivations, inspired  from the field of 
region-based theory of space  and parallel to the arguments
presented here in the introduction, suggesting that it is natural to deal
with the join operation and possibly excluding meets and complements.  
The paper \cite{DIM}  analyzes aspects of 
pointfree topology which can be expressed in terms
of a particular class of join semilattices; the difference with the present note
is that in \cite{DIM} semilattices have no additional structure.
\arxiv{For complete lattices,  \cite[Sects.\ 2, 3]{Now}  deals with  binary relations
defined by condition \eqref{1} 
 associated
 to operators satisfying 
properties weaker  than closure. }

\section{Embedding theorems} \labbel{embsec} 

\subsection{Principal specialization semilattices} \labbel{prin}

In our motivating example  from
Definition \ref{asso} the relation
$ a \sqsubseteq b$ is defined by $ a \subseteq   \K b$,
where $ \K $ is the topological closure of some topological space.
In particular, in the motivating example, $ \K b$ is the largest element of the set
$S_b=\{ a \in S \mid  a \sqsubseteq b \}$.
On the other hand, in a general specialization poset,
as we have defined it,
 $S_b$ might not even have a  (possibly infinitary) 
join. \arxiv{See Example \ref{ord}(a) below. }
Even when $S_b$ has a join, say, $s$,
it is not necessarily the case that 
$s \sqsubseteq b$.
\arxiv{See Example \ref{ord}(b), 
or consider inclusion mod finite on an
infinite set. }
In the following lemma we show that if, for some $b \in S$,
the set  $S_b$
has a maximum, call it $ \K b$, 
then $ \K b$ actually corresponds to some form of 
``closure''. The lemma works for a single $b$
(possibly, two elements $a$ and  $b$) and we do not need
the assumption that $S_b$ has a maximum for \emph{every}
$b \in S$. 
\arxiv{If the latter property actually holds 
for every $b \in S$, we get a very special class 
(widely known, in an equivalent formulation)
 of specialization posets we are going to mention soon. }
In some special situations the following lemma is
known; for example, the equivalence of (ii) and (iv)
has been used in order to provide an equivalent definition 
for a closure operation \cite[Proposition 3.1]{E}. 
 The lemma is probably new at the present level
of generality.

\begin{lemma} \labbel{lemlemlem}
Suppose that $\mathbf S$ is a specialization poset,
$b \in S$ and the set  
$\{ a \in S \mid  a \sqsubseteq b \}$
has a maximum (in the sense of $\leq$). 
Call this maximum $ \K b$. Then, for every $a \in S$,
the following conditions are equivalent.
 \begin{enumerate}[(i)]
   \item
$ a \sqsubseteq b $;
\item
$ a \leq  \K b$;
\item
$ a \sqsubseteq  \K b $.
   \end{enumerate}   
 In particular, $ \K b$
is also the maximum 
of $\{ a \in S \mid  a \sqsubseteq  \K b \}$.
Suggestively, if $ \K b$ exists, then  $ \K  \K b $ exists, too, and 
they are equal. 

Suppose further that $a \in S$ and the set 
$\{ c \in S \mid  c \sqsubseteq a \}$
has a maximum, call it $ \K a$. 
Then conditions (i) - (iii) above are also equivalent to
 \begin{enumerate}[(i)]
\setcounter{enumi}{3}
 \item
$  \K a \leq  \K b$;
\item
$  \K a \sqsubseteq  \K b $.
   \end{enumerate}   
 \end{lemma} 

\begin{proof}
The fact that  (i) implies (ii)
is just a restatement of the definition of $ \K b$.
From (ii) we immediately get (iii),
because of \eqref{s1}.
By the definition of $ \K b$
we get $ \K b \sqsubseteq b$,
hence from (iii) and \eqref{s2}
we obtain (i).
Under the additional assumptions,
$ \K a \sqsubseteq a$, hence (iii)
implies (v), by \eqref{s2}.
Moreover,  (iv) and (v) are equivalent,
by applying the equivalence of (ii) and (iii)
with $ \K a$ in place of $a$.  
Since $a \leq  \K a$, by \eqref{s4},
we get that (iv) implies (ii).  
\end{proof}    

\arxiv{The above considerations and 
Lemma \ref{lemlemlem}  suggest the following well-known definition. }

\begin{definition} \labbel{supcomp}   
(a) A specialization poset (semilattice)
$\mathbf S$ 
is \emph{principal} 
if the preorder $ \sqsubseteq $
is \emph{a principal quasi-order for $(S, \leq)$}
\cite[p. 193]{E}.
This means that,
 for every   $b \in S$,
the set 
$S_b = \{ a \in S \mid  a \sqsubseteq b \}$
has a maximum relative to $\leq$.
In other words, $\mathbf S$
is principal if and only if every principal $ \sqsubseteq $-ideal is also 
$\leq$-principal.    
 
(b) In the case of a specialization semilattice, 
$S_b$ is (finitely) upward directed, by \eqref{s3}.
In particular, every finite 
specialization semilattice is principal.
\arxiv{However, as we mentioned, 
when $S_b$ is infinite
it is not necessarily the case 
that $S_b$ has a  
join. Even when some join $s$
exists,
it is not necessarily the case that 
$s \sqsubseteq b$.
In the definition of a principal quasi-order
we not only  require that such an $s$ exists, but 
we require that actually
$s \sqsubseteq b$.  }
On the other hand, a finite specialization poset
is not necessarily principal.  Example \ref{diamond}
furnishes a counterexample: in that example
$S_0 = \{ x\in S \mid  x \sqsubseteq 0 \}$
has no maximum.

\arxiv{(b$'$) In particular, if $\mathbf S$ 
is a specialization semilattice which is complete as a 
semilattice, then  $\mathbf S$ is principal 
if and only if $\mathbf S$ satisfies the following infinitary version
of condition \eqref{s3}, for every
family $( a_i) _{i \in I} $ of elements of $S$.
\begin{equation}
\labbel{s3oo}    \tag{S3$_{\infty}$}
 \text{If $a_i \sqsubseteq b  $, for every $i \in I$,
then  
 $\bigvee_ {i \in I}  a_i \sqsubseteq b $.}
\end{equation}

Of course, we do not need the full assumption that $\mathbf S$  is complete: it is
enough to assume that all the joins as given by \eqref{s3oo}  exist.

Principal specialization complete lattices are 
\emph{complete lattices with a symmetric consequence relation} 
in the terminology from \cite[3.1]{GT} (considering the converse
of the consequence relation).  }

(c) If $\mathbf S$ is a principal specialization poset (semilattice),
then, for every $b \in S$, 
let us denote by $ \K b$ the maximum of $S_b$.
Then the  function 
assigning $ \K b$ to $b$  is a closure operation
\cite[Proposition 3.9]{E}.

We say that a principal specialization semilattice is
\emph{additive} if $ \K (a \vee b) =  \K a \vee  \K b$,  
for every $a,b \in S$. 
 \end{definition}

\arxiv{
\begin{remarkx} \labbel{sc}
(a) The specialization  semilattice (poset)
associated to a closure space,
as in Definition \ref{asso}, is principal.
The specialization semilattice associated to some topological space 
is also additive. Indeed, a closure space is a topological space
if and only if its associated specialization  semilattice 
is additive in the above terminology, and moreover 
$ \K  \emptyset = \emptyset $.
 In particular, not every principal specialization semilattice 
is additive.

(b) If $\mathbf S$ is a principal specialization poset (semilattice),
then the  function 
which assigns to $b$ the maximum $ \K b$
of $S_b = \{ a \in S \mid  a \sqsubseteq b \}$ is a closure operation,
by Lemma \ref{lemlemlem}.
Recall Definition \ref{closop}. 
Conversely, if $ \K $ is a closure operation
on some poset (semilattice), then, as noted in Remark \ref{si!},
we get a principal specialization poset (semilattice) 
by letting 
$ a \sqsubseteq b$ if 
$a \leq  \K b$. Notice that 
the motivating example from  
Section \ref{motiv} 
 is a special case of the above construction.

The above constructions provide a bijective correspondence
between
closure operations and principal specializations on the same poset. 
See
\cite[Proposition 3.9]{E}.
Let us point out (again)
that the notions of homomorphism are distinct in the
two settings.
 \end{remarkx}

Remark \ref{sc}(b) can be used  to show
that specializations provide still another equivalent formulation
for the notion of a topology, possibly a folklore result
in some form or another. 

\begin{observationx} \labbel{obs}
Fix some set $X$.
 The correspondence assigning 
to some closure\footnote{Here, of course, by \emph{closure}
we mean the family of the closed sets of a closure 
space.}  (topology) $\tau$ on $X$ 
the 
specialization semilattice 
$\mathbf S(X, \tau )$ from Definition \ref{asso} 
is a bijective correspondence from the set of all
closures  (topologies) on $X$
to the set of all the principal (principal and additive)
 specialization semilattices of the form
 $( \mathcal P(X), \cup, {\sqsubseteq} )$
(and such that $ a \sqsubseteq \emptyset $ implies $a= \emptyset $).  
 
In particular, by Proposition \ref{propcont},
 the category of topological spaces 
with continuous functions is isomorphic  
as a concrete category over \textbf{Set} \cite[Remark 5.12]{joy} 
to the
category of 
principal and additive
 specialization semilattices which can be realized as
 $( \mathcal P(X), \cup, {\sqsubseteq} )$, for some set $X$,
and such that $ a \sqsubseteq \emptyset $ implies $a= \emptyset $,
with homomorphisms.
 \end{observationx}

\begin{remarkx} \labbel{scbis}
(a)
Notice that a principal specialization semilattice is additive if and only if
\begin{equation}\labbel{kuk}     
\text{$ \K (e \vee d) = e \vee d$, for every $e,d$ such that 
$ \K e=e$ and  $ \K d=d$.}
 \end{equation}
Necessity is obvious.
In the other direction,
for every $a$ and $b$,  
take $e= \K a$ and $d= \K b$.
We have   $ \K e=e$ and  $ \K d=d$ by Lemma \ref{lemlemlem}.
Then \eqref{kuk} provides 
 $ \K ( \K a \vee  \K b) =  \K a \vee  \K b$,
hence
 $ \K (a \vee b) \leq   \K ( \K a \vee  \K b) =  \K a \vee  \K b$,
the other inequality being obvious,
since $ \K $ is isotone, due to \eqref{s5}.

(b) The condition that some specialization semilattice (poset)
is principal can be expressed by a first-order sentence. Indeed,
a specialization semilattice (poset) $\mathbf S$ is principal if and only if 
$\mathbf S$ satisfies
\begin{equation}\labbel{p}    
\forall b \exists c \forall a ( a \sqsubseteq b \Leftrightarrow a  \leq c).
  \end{equation}
Thus, \emph{provided some specialization semilattice $\mathbf S$ is principal},
by (a)  we can express the property that $\mathbf S$ is additive as 
\begin{multline}\labbel{pp}    
\forall  cd (
\forall a  ( a \sqsubseteq c \Leftrightarrow a  \leq c)
\ \& \ 
\forall a  ( a \sqsubseteq d \Leftrightarrow a  \leq d)
\Rightarrow 
\\
\forall a  ( a \sqsubseteq c \vee d \Leftrightarrow a   \leq c \vee  d)
).
  \end{multline}

However, the sentences \eqref{p} and
\eqref{pp} are  quite complex; here we are interested in simpler sentences,
mainly universal sentences.  
Of course, if we add further symbols
 to our language, the statement that a specialization semilattice 
is principal
can be expressed in a simpler way.
If we add the unary function $ \K $,
then
\eqref{p} can be expressed as 
$\forall ab  ( a \sqsubseteq b \Leftrightarrow a  \leq  \K b)$.
However, $ \K $ is not definable over every
specialization semilattice (poset);
in fact, we have $ \K $ exactly when the specialization semilattice 
(poset) is principal.
 If we abbreviate 
$\forall a  ( a \sqsubseteq c \Leftrightarrow a  \leq c)$
as $ \CC (c)$, then 
\eqref{pp} simplifies to   
$\forall  cd ( \CC (c) \ \& \  \CC (d) \Rightarrow \CC (c \vee d))$. 

However, as commented above,
changing the language changes the notion
of homomorphism. In particular, if we add, say, a unary predicate $ \CC $ 
as defined above, the analogue of Proposition \ref{propcont}
fails, since not every continuous function is closed.
\end{remarkx}   
 }

\subsection{Embedding into principal specialization semilattices} \labbel{emb} 

\arxiv{
Homomorphic images of specialization semilattices 
are not necessarily specialization semilattice themselves.
 See Example \ref{exquot}
for an explicit  counterexample. 
Technically, this can be hinted from  the fact that 
\eqref{s1} - \eqref{s3} are not positive formulas, 
and  shows that 
the theory of specialization semilattices cannot be axiomatized 
by positive sentences, by \cite[Theorem 3.2.4]{CK}.
Hence   Lemma \ref{quot} below will be useful.
We first recall a classical algebraic definition. }

\begin{definition} \labbel{cong}    
If $\mathbf S =(S, \vee )$
is a semilattice, a binary relation
$ \sim$ is a \emph{congruence}
on $\mathbf S$ 
if $ \sim$ is an equivalence relation 
on $S$ 
and furthermore
$ \sim$ respects $\vee$, that is,
$ a \sim b$ implies $ a \vee c \sim b \vee c $,
for every $a,b,c \in S$.   
This is a special case of a more general
algebraic notion \cite[Section 1.5]{B}. 

If $\mathbf S$ has some further structure,
we will say that $ \sim$ is \emph{a congruence 
 for the semilattice reduct} if the above conditions hold.
\end{definition}

\begin{lemma} \labbel{quot}
Suppose that $\mathbf S =(S, \vee, {\sqsubseteq} )$
is a specialization semilattice and
$ \sim $ is an equivalence relation on $S$ such that
$ \sim$ is a congruence for the semilattice reduct
 and moreover 
\begin{equation}\labbel{seq}
\text{ for every $a,b \in S$,
if $ a \sim b$,  then $ a \sqsubseteq b $ 
and $ b \sqsubseteq a$.} 
   \end{equation}

Then $\undertilde{\mathbf  S} =(\undertilde{S} ,
 \mathbin {\undertilde{\vee}} , \mathrel {\undertilde{\sqsubset}}  )$
is a specialization semilattice, 
where $\undertilde{S} $  is the set of the $ \sim$-equivalence classes,
$\mathbin {\undertilde{\vee}} $ is the standard quotient operation
and,  for all $a,b \in S$,  
we let $ \undertilde{a}  \mathrel {\undertilde{\sqsubset}} \undertilde{b}   $
if $a \sqsubseteq b$.
Here we have written, say, $\undertilde{S}$, $\undertilde{a}$, \dots\ 
in place of   $S/{\sim}$, $a/{\sim}$ or $[a]$,\dots\ 
in order to improve readability.

Moreover, the projection map 
$\pi$ which sends $a$ to $\pi(a)= \undertilde{a}$ is a homomorphism
of specialization semilattices.  
 \end{lemma} 

\begin{proof}
First, we 
show that
$\mathrel {\undertilde{\sqsubset}}$  is well-defined.
In fact, if $a \sim a_1$, 
$b \sim b_1$ and $ a \sqsubseteq b$,
then by \eqref{seq}  we have
$ a_1 \sqsubseteq a$ and 
$ b \sqsubseteq b_1$, hence 
 $ a_1 \sqsubseteq b_1$ by \eqref{s2}. 

By standard algebraic arguments \cite{B} 
 $\mathbin {\undertilde{\vee}} $ is a semilattice operation on $\undertilde{S} $,
since we have assumed that $ \sim$ is a  congruence
for the semilattice reduct; moreover, the projection
is a semilattice homomorphism.

If $ \undertilde{a} \mathrel{\undertilde{<}} \undertilde{b}$
in  $ \undertilde{\mathbf  S} $, 
then, by the above paragraph, 
$ (a \vee  b)/{\sim} = 
\undertilde{a} \mathbin {\undertilde{\vee}} \undertilde{b} = \undertilde{b} $,
that is,
$ a \vee  b \sim b $.
By \eqref{seq}
 $ a \vee  b \sqsubseteq  b $,
hence 
$ a  \sqsubseteq  b $ by \eqref{s6},
thus
$ \undertilde{a}  \mathrel {\undertilde{\sqsubset}} \undertilde{b}   $,
by the definition of $\mathrel {\undertilde{\sqsubset}}$.
We have proved that
\eqref{s1} holds 
in  $\undertilde{S} $.

If $ \undertilde{a} \mathrel {\undertilde{\sqsubset}} \undertilde{b}$
and 
 $ \undertilde{b} \mathrel {\undertilde{\sqsubset}} \undertilde{c}$
then 
$a \sqsubseteq b$ and  $b \sqsubseteq c$,
since the definition of $ \mathrel {\undertilde{\sqsubset}} $
does not depend on the representatives.
Hence $a \sqsubseteq c$ by \eqref{s2}
and   $ \undertilde{a} \mathrel {\undertilde{\sqsubset}} \undertilde{c}$
by the definition of $ \mathrel {\undertilde{\sqsubset}} $,
so that \eqref{s2}  holds in $\undertilde{\mathbf  S}$ as well. 

The proof of \eqref{s3} is similar,
using the already mentioned fact that
$ (a \vee  a_1)/{\sim} = 
\undertilde{a} \mathbin {\undertilde{\vee}} \undertilde{a_1}  $.
The last statement is trivial. 
\end{proof}    

\arxiv{
The assumption that $ \sim$ satisfies Condition \eqref{seq} 
is necessary in Lemma \ref{quot}. See Example \ref{exquot}
below.  }

\begin{theorem} \labbel{propcomp}
Every specialization semilattice
can be embedded into  a 
  principal additive
 specialization semilattice.
 \end{theorem}  

\begin{proof}
Suppose that $\mathbf S =(S, \vee_S, \sqsubseteq_S )$
is a specialization semilattice.  
Let $T= \{ 0,1 \}$  and  $\mathbf T= (T, \max, \sqsubseteq_T ) $ 
be the specialization semilattice 
such that $ x \sqsubseteq_T y$ for every $x,y \in T$.
Let $\mathbf S \times \mathbf T$ be defined in the natural way
on $S \times T$, by taking the standard  semilattice product  
  and letting $ a \sqsubseteq _{S \times T} b $
hold if both
$ a_1 \sqsubseteq_S b_1$ and 
$ a_2 \sqsubseteq_T b_2$,
for all $a= (a_1, a_2)$ and $b= (b_1, b_2)$
in $S \times T$
(of course, in the special case at hand, 
$ a \sqsubseteq _{S \times T} b $
holds if and only if 
$ a_1 \sqsubseteq_S b_1$ holds).
By Remark \ref{conseq}(d),  
$ \mathbf S \times \mathbf T$ is a specialization semilattice.

The function $\iota$ defined by
 $\iota(s) = (s,0)$
is a homomorphism from $\mathbf S$ to 
$\mathbf S \times \mathbf T$.
This fact easily follows from  the observation
that $\{ 0 \}$ is the universe for a submodel of 
$\mathbf T$ and, of course, using 
$0 \sqsubseteq_T 0$. 

Next, consider the equivalence relation
$\sim$ on  
$  S \times  T$
defined by 
$(a_1, a_2) \sim (b_1, b_2)$ 
if at least one of the following conditions hold
  \begin{enumerate}[(1)]
    \item  
 $a_1=b_1$ and  $a_2=b_2$, or 
\item
$a_2 = b_2=1$ and both
$ a_1 \sqsubseteq_S b_1 $ and $ b_1 \sqsubseteq_S a_1 $.
  \end{enumerate}

Due to \eqref{s2}, $ \sim$ is transitive,
and then an equivalence relation by \eqref{s4}. 

We now check that $ \sim$ is a congruence for the semilattice
reduct. Suppose that
$(a_1, 1) \sim (b_1, 1)$, as witnessed by 
the conditions
$ a_1 \sqsubseteq_S b_1 $ and  $ b_1 \sqsubseteq_S a_1 $.
Then, for every $(c_1, c_2) \in S \times T$, we have 
$(a_1, 1) \vee (c_1, c_2) =
(a_1\vee_S c_1, 1)$ and 
$(b_1, 1) \vee (c_1, c_2) =
(b_1\vee_S c_1, 1)$.
From $ a_1 \sqsubseteq_S b_1 $
and  \eqref{s9}
we get $a_1\vee_S c_1 \sqsubseteq _S b_1\vee_S c_1$ 
and, symmetrically,
$ b_1\vee_S c_1 \sqsubseteq _S a_1\vee_S c_1$, 
thus  $(a_1\vee_S c_1, 1) \sim (b_1\vee_S c_1, 1)$,
that is, 
$(a_1, 1) \vee (c_1, c_2) \sim (b_1, 1) \vee (c_1, c_2)$.
The other case is  trivial.
We have showed that $ \sim$ is a congruence for the semilattice reduct.

Finally, $ \sim$ satisfies Condition \eqref{seq} 
by construction,
hence we can apply Lemma \ref{quot}
in order to get a specialization semilattice  
$\mathbf U = (\mathbf S \times \mathbf T) / {\sim} $
and a homomorphism
$\pi : \mathbf S \times \mathbf T \to \mathbf U$.
The composition 
$\kappa = \pi \circ \iota$,
being a composition of two homomorphisms,
is a homomorphism from $\mathbf S$ to  $\mathbf U$.  

In order to keep the number of 
subscripts to a minimum, let us write 
$[a_1, a_2]$
for $ \pi (a_1, a_2)$, in place of 
$(a_1 \undertilde{,}\,  a_2 )  $ or $(a_1, a_2) /{\sim}$.
The homomorphism $\kappa$ 
is injective, since
if $a \neq b \in S$,
then 
$ \kappa (a) = [a, 0]$  
and  
$ \kappa (b) = [b, 0]$,
thus
$[a, 0] \neq  [b, 0]$,
since no two pairs with  second component 
$0$ are identified by $ \sim$.
By the definition of $ \mathrel {\undertilde{\sqsubset}} $
if  $[a, 0] \mathrel {\undertilde{\sqsubset}} [b, 0]$,
then  
$(a, 0) \sqsubseteq_{S \times T}  (b, 0) $
and then
$a \sqsubseteq_{S}  b $, by the definition 
of $\sqsubseteq_{S \times T}$.
Thus $\kappa$ is an embedding of 
$\mathbf S$  into $\mathbf U$.

It remains to show that 
$\mathbf U$ is principal and additive.
Taken any element 
$[a_1, a_2]$ 
of $\mathbf U$,
we see that
 $[a_1, 1] \mathrel {\undertilde{\sqsubset}} [a_1, a_2]$,
by the definition of $ \mathrel {\undertilde{\sqsubset}}$ and since 
$(a_1, 1) \sqsubseteq_{S \times T} (a_1, a_2) $,
by \eqref{s4} and the definitions of
$ \sqsubseteq_T $
and   of $ \sqsubseteq _{S \times T}$. 
We claim that
 $[a_1, 1] $ is the largest element
  $[c, d] \in U $ such that 
 $[c, d] \mathrel {\undertilde{\sqsubset}} [a_1, a_2]$.
Indeed, if 
 $[c, d] \mathrel {\undertilde{\sqsubset}} [a_1, a_2]$,
then
 $(c, d)  \sqsubseteq_{S \times T} (a_1, a_2) $, hence
 $c  \sqsubseteq_{S } a_1 $,
by the definitions of $\mathrel {\undertilde{\sqsubset}}$
 and $\sqsubseteq_{S \times T}$.
It follows that
 $(c, 1)  \sqsubseteq_{S \times T} (a_1, 1) $,
hence 
 $(c, 1) \vee_{S \times T} (a_1, 1)
  \sqsubseteq_{S \times T} (a_1, 1) $,
by \eqref{s8}.
By \eqref{s1},
   $ (a_1, 1)  
  \sqsubseteq_{S \times T} (c, 1) \vee_{S \times T} (a_1, 1)$,
hence
 $(c, 1) \vee_{S \times T} (a_1, 1)
  \sim (a_1, 1) $,
by the definition of $ \sim$.
 Thus in $\mathbf U$ 
\begin{equation*}    
 [c, 1] \mathbin {\undertilde{\vee}} [a_1, 1]
=\pi (c, 1)\mathbin {\undertilde{\vee}} \pi (a_1, 1)=
\pi ((c, 1) \vee_{S \times T} (a_1, 1))=
 [a_1, 1] ,
 \end{equation*}
 since $\pi$
is a semilattice homomorphism.
It follows that in $\mathbf U$ 
$ [c, 1] \leq [a_1, 1]$.
Obviously, 
$[c, d] \leq [c, 1]$,
hence 
$ [c, d]\leq [a_1, 1]$.

 We have proved that 
 $[a_1, 1] $ is the maximum among those  elements
  $[c, d] $ such that 
 $[c, d] \mathrel {\undertilde{\sqsubset}} [a_1, a_2]$.
In the notation from
Definition \ref{supcomp}(c),
 $ \K [a_1, a_2] = [a_1, 1]$. 
Since $[a_1, a_2]$ is arbitrary,
the above procedure 
applies to every element of $U$,
thus $\mathbf U$ 
is principal.
 Moreover $\mathbf U$ is additive, since
\begin{align*}
\K ([a_1, a_2] \mathbin {\undertilde{\vee}} [b_1,b_2]) &=
 \K ( \pi(a_1, a_2) \mathbin {\undertilde{\vee}} \pi(b_1,b_2))=
 \K (\pi((a_1, a_2) \vee_{S \times T} (b_1, b_2))
\\
&= \K \pi(a_1 \vee_S b_1, a_2 \vee_T b_2)=
 \K [a_1 \vee_S b_1, a_2 \vee_T b_2]
\\
& = [a_1 \vee_S b_1, 1] =
[a_1 , 1] \mathbin {\undertilde{\vee}}   [ b_1, 1] =
 \K [a_1 , a_2]  \mathbin {\undertilde{\vee}}  \K [b_1 , b_2].  \qedhere
  \end{align*}    
 \end{proof} 

\arxiv{
\begin{remarkx} \labbel{prac}   
(a) In practice, the extension $\mathbf U$ 
 in the above proof is obtained as follows.
Consider the equivalence relation $ \Theta $
on $S$ defined by 
$a \mathrel { \Theta }  b$
if $ a \sqsubseteq_S b $ and  $b \sqsubseteq_S a$.   
Since $\Theta$ satisfies \eqref{seq}, by \eqref{s9}, the quotient $S/{\Theta}$ 
naturally becomes a specialization semilattice,
call it $\mathbf V$.
Then $\mathbf U$ can be thought of as the
``disjoint union'' of  
$\mathbf S$ and $\mathbf V$,
obtained by declaring each class in 
$V$ to be larger than (actually, the closure of) each
member of the class.
While the proof of Theorem \ref{propcomp}
can be performed along the above lines,
the proof becomes harder,
since there are a lot of details to be checked by hand. 
Most of these details automatically 
follow from the canonical structures
on products and quotients, as presented in 
the  given proof of Theorem \ref{propcomp}.

(b)
The embedding $\kappa$ in the proof of \ref{propcomp}
preserves existing, possibly infinitary,  meets 
of nonempty subsets of $\mathbf S$.
Indeed, in $\mathbf U$ it never happens that
$[a,1] \leq [b,0]$, hence existing meets are computed as in     
$\mathbf S$. 

On the other hand, $\kappa$  never
preserves existing 
closures, since a new closure in $\mathbf U$
is always added, even for 
elements already having a closure
in S. See \cite[Sect. 4]{mttna} 
for a construction preserving any
prescribed set of existing closures.

(c)
In general, there is   not a smallest extension
of  
$\mathbf S$ satisfying the conclusions of Theorem \ref{propcomp}.
See Example \ref{star}. On the other hand, the construction in the 
proof of Theorem \ref{propcomp} provides a canonical  ``free''
extension in the class of closure spaces satisfying 
$a \vee Kb=K(a \vee b)$. See Case (C5) 
in the final section of \cite{mttlib}. 
\end{remarkx} 
 }  

\subsection{The universal theory of topological spaces
(in the language of specialization semilattices)} \labbel{univcs}

Recall from
Definition \ref{asso}
that to a topological space $X$ 
 we have associated 
the specialization semilattice 
$\mathbf S(X) = (\mathcal P(X), \cup , \sqsubseteq  )$
and 
the specialization poset 
$\mathbf P(X) = (\mathcal P(X), \subseteq  , \sqsubseteq  )$,
where $a \sqsubseteq b$ if 
$ a \subseteq  \Ku b$, $ \Ku $ being the closure induced by 
the topology on $X$. 
We say that a specialization semilattice, resp., poset
is \emph{topological} if it is isomorphic to
$\mathbf S(X )$,  
resp., to $\mathbf P(X )$,
for some topological space $X$. 
Similarly, the \emph{closure semilattice
associated} to some closure or topological space $X$ 
 is   $(\mathcal P(X), \cup  , \K )$,
where $\K$ is the closure of $X$.

\begin{proposition} \labbel{clemb}
Every (additive) closure semilattice 
can be embedded into the 
closure semilattice associated to a closure (topological) space.
 \end{proposition}  

\begin{proof}
Let $\mathbf S = (S, \vee, \K _S )$
be a closure semilattice.
First notice that we can assume
that $\mathbf S$ has a minimum element $0$
such that $\K _S 0 =0$. If not,     
$\mathbf S$ can be embedded into a closure semilattice 
$\mathbf S_0$ with such an element; just add 
to $\mathbf S$
a new $\vee$-neutral element
$0$ and set $\K 0 =0$.
Then the inclusion map is an embedding
of $\mathbf S$ into $\mathbf S_0$.
If $\mathbf S$ is additive, then $\mathbf S_0$
is additive, as well. 

Let $\varphi: S \to \mathcal P(S)$
be the function defined by  
$ \varphi (a) = \nup a = 
\{ b \in S \mid a \nleq b \}  $. 
 We will define a closure on  
$S$  in such a way that  
 $\varphi$  is an embedding
from $\mathbf S$ to the associated closure semilattice.
First, $\varphi$  is  an injective  semilattice homomorphism,
hence an embedding
from  $(S, \vee)$ to $(\mathcal P(S), \cup)$, no matter how
the closure on $S$  is defined. 

Let $S_1 \subseteq \mathcal P(S)$ 
be the image of $S$ under $\varphi$  and
define a function $\K_1: S_1 \to \mathcal P(S)$  
as follows.  
If $x= \varphi (a)$,  let $\K_1 x= \varphi ( \K _S a)$.
$\K_1$ is well-defined since $\varphi$  is injective.
Since $\nup 0 = \emptyset $, then (i) $ \emptyset \in S_1$
and $\K_1 \emptyset = \emptyset $, since $\K _S 0 = 0$.    
By construction, (ii) if $x \in S_1$, then $\K _1 x \in S_1$
and   $x \subseteq \K _1 x = \K _1 \K _1 x$, since
 $\K _S$ is a closure on $\mathbf S$ and $\nup$ is 
order preserving. If $\K _S$ is additive, then
(iii) if  $x,y \in S_1$, then $ x \cup y  \in S_1$
and   $\K _1 (x \cup y) = \K _1 x \cup  \K _1 y$,
by additivity of $\K _S$ 
and since $\nup$ is a semilattice homomorphism.
  
The above conditions (i) - (iii) are exactly the sufficient conditions
in \cite[Lemma 2.3]{MT}  for extending $\K_1$ to an  additive closure operation 
$\K$ on the whole of  $\mathcal P(S)$, which is
a completely additive Boolean algebra.
Since $\mathcal P(S)$ is the field of all the subsets of $S$,
then $\K$ is a topological closure on $S$; moreover, $\varphi$  is a
$\K$-embedding by construction. In conclusion, $\varphi$   is an embedding
of closure semilattices. The above argument treats the additive case;
the nonadditive case  is slightly simpler and  uses \cite[Lemma 8.1]{Sc}.
See also \cite[Proposition 3]{ecca} and \cite[Lemma 3.1]{sapimpap}.
 \end{proof}

\begin{theorem} \labbel{sptop}
Every specialization  semilattice
can be embedded into a topological
specialization semilattice.  
 \end{theorem}

\begin{proof}
By Theorem \ref{propcomp}
every specialization semilattice $\mathbf S$ can be embedded into an additive
principal specialization semilattice. By adding the corresponding closure operation
as in Definition \ref{supcomp}(c),
we get an additive closure semilattice, which, by Proposition \ref{clemb},
 can be embedded into the
closure semilattice associated to some topological space.
Since, as noticed in a comment in Definition \ref{closop},
  an embedding of closure semilattices is an embedding
of the associated specialization semilattices, then, by composing the
two embeddings,
we get an embedding of $\mathbf S$ into a topological
specialization semilattice.
\end{proof}   

\arxiv{
Notice that, by Remark \ref{sc}(a), 
Theorem \ref{sptop}
is formally stronger than  Theorem \ref{propcomp};
however, the  proof  of \ref{sptop} makes use of  \ref{propcomp}.

\begin{remarkx} \labbel{nomeet}    
In contrast with Remark \ref{prac}(b),
the embedding $\varphi$  in the  proof of Theorem \ref{sptop} 
does not generally preserve existing meets.

In fact, preservation of meets cannot be accomplished,
in general.
Just consider some lattice
$\mathbf L$ which is not distributive \cite{Bly} and  let $ \K $ be the identity
function on $L$. The associated specialization semilattice 
is principal and additive,
but it cannot be embedded into the specialization semilattice 
associated  to some topological space in such a way 
that meets are preserved, since every sublattice of a
 distributive lattice is distributive \cite{Bly}. 
\end{remarkx}
 }

By Remark \ref{conseq}(d)
the class of specialization semilattices 
(posets) is axiomatized by a first-order universal theory,
which we will call the \emph{theory of specialization semilattices
(posets)}.  
Recall the definitions of $\mathbf S(X)$
and  of $\mathbf S(P)$
from Definitions \ref{asso} and \ref{closop}.

\begin{corollary} \labbel{corun}
For every universal first-order sentence $\sigma$
in the language of specialization semilattices,
the following conditions are equivalent. 
\begin{enumerate}[(1)]   
\item   
The sentence $\sigma$  is true
in all the models of the form
$\mathbf S(X)$, where $X$ varies among topological spaces.
\item
The sentence $\sigma$  is true
in all the models of the form
$\mathbf S(X)$,  where $X$ varies among closure spaces.
\item
The sentence $\sigma$  is true
in all the models of the form
$\mathbf S(P)$,  where $P$ varies among closure semilattices.
\item
The sentence $\sigma$
 is a logical consequence
of the theory of specialization semilattices.
\end{enumerate}
 \end{corollary}

 \begin{proof}
(1) $\Rightarrow $  (4) If $\sigma$  
is valid in all the models  $\mathbf S(X)$,
with $X$ a topological space, then
$\sigma$  is valid in all the specialization semilattices,
by Theorem \ref{sptop}
and since universal sentences are preserved under
taking submodels. 

(4)\hspace{1 pt}$\Rightarrow $\hspace{1 pt}(3) follows from 
Remark \ref{si!};
(3)\hspace{1 pt}$\Rightarrow $\hspace{1 pt}(2)
and (2)\hspace{1 pt}$\Rightarrow $\hspace{1 pt}(1) are  obvious.
 \end{proof}  

For short, the universal theory of specialization semilattices
is the universal theory of 
structures associated to topological spaces 
in the sense of Definition \ref{asso}. 

\begin{remark} \labbel{subt} 
The equivalences of (1) - (3)
in Corollary \ref{corun}
assert  that \emph{in the language of specialization semilattices}
the universal sentences valid in all $\mathbf S(X)$
are the same, no matter whether we let 
$X$ vary among topological spaces,  closure spaces,
or even closure semilattices.

Thus Corollary \ref{corun} might be interpreted 
in the sense that the language $\{  \vee, \sqsubseteq \}$
lacks  expressive power since, as far as universal sentences
are considered, it does not distinguish between 
subreducts of topological spaces
and  subreducts of closure spaces. However, the issue is subtler,
as we are going to explain.

If $X$ is a topological space, or just a closure space,
 with closure operator $ \K $,
consider the following model
 $\mathbf M(X) = ( \mathcal P(X), \cup, R)$,
where $R$ is the ternary relation defined by
\begin{equation*}   
\text{  $R(a,b,d)$ \quad if  \quad  
$a  \subseteq   \K b \cup  \K d$. }
 \end{equation*}
As in Proposition \ref{propcont},
 some function   $\varphi$  is  continuous   
between two  topological (or closure) spaces $X$ and $Y$
 if and only if $\varphi^\sharr$
is a homomorphism between the corresponding models
  $\mathbf M(X)$ and $\mathbf M(Y)$.

Notice that 
$a \sqsubseteq b$
from Definition \ref{asso} 
is interpretable as
$R(a,b,b)$.
Since in a topological space 
we have
$ \K (b \cup d) =  \K b \cup  \K d$, then, if $X$ 
is a topological space, the following 
equivalences hold  in $\mathbf M(X)$:
\begin{equation}\labbel{buh}  
R(a,b,d)   \Leftrightarrow  R(a,b \vee d ,b \vee d)
\Leftrightarrow a \sqsubseteq b \vee d,
 \end{equation}    
for all $a,b,d \subseteq X$.
On the other hand, the first equivalence in \eqref{buh}
does not necessarily hold in 
$\mathbf M(X)$, when $X$ is just assumed to be a closure
space. 
Simply take $b$ and  $d$
to be two closed whose union $ b \cup d$ is not closed 
 and let $a=  \K (b \cup d)$, thus $a \supsetneq  b \cup d$.
Then   
$  R(a,b \vee d ,b \vee d) $ holds in $\mathbf M(X)$, but 
$R(a,b,d) $ fails, hence \eqref{buh} fails, as well.  
Actually, the argument shows that $ \K $ is additive if and only if
the first equivalence in \eqref{buh} holds.

Thus in the language $\{ \vee, R\}$
the universal theories of topological spaces 
and of closure spaces are distinct.
However, \emph{when we deal with topological spaces},
the relation $R$ can be defined in terms of $ \sqsubseteq $,
by the equivalence of the outer conditions in \eqref{buh},
hence it is sufficient to deal with $ \sqsubseteq $. 
 In other words, the language $\{  \vee, \sqsubseteq \}$
does distinguish topological spaces and closure spaces under
a definitional  expansion which,
 \emph{in the case of topological spaces}, is unessential.
In practice,  if we think
of $ a \sqsubseteq b \vee d$ as $a \subseteq  \K b \cup  \K d$,
this interpretation is sound for topological spaces but not for 
closure spaces.  

Hence  Corollary \ref{corun} only shows that 
 the language $\{  \vee, \sqsubseteq \}$ lacks
expressive power \emph{when the intended examples
are closure spaces}. Under an unessential definitional expansion,
that is, adding $R$, we do distinguish topological spaces 
from non-topological closure spaces.
The above remark is also supported by the observation that 
the construction of ``universal'' closure semilattices
is much simpler in the additive case (Theorem \ref{propcomp} or \cite{mttlib})
with respect to the nonadditive case \cite{mttna}.   
\arxiv{

We need $\vee$ in order to perform
 the above arguments. It is not clear how expressive the language
$\{  \leq , \sqsubseteq \}$ is, in the context of topological spaces. }
\end{remark}

\arxiv{The above remarks suggest the following problems.

\begin{problemx} \labbel{ur}
  \begin{enumerate}[(1)]    
\item  
 With
the above interpretations,
provide axioms for the  universal theories 
of topological spaces and of closure  spaces in the  language
 $\{ \leq, R\}$, as well as of closure
spaces in the language  $\{ \vee, R\}$.

\item
Does any qualitative difference arise
if we add  $n+1$-ary relations $R_n$ whose intended  interpretation is 
$R_n(a,b_1,b_2, \dots, b_n)$ if 
$a  \leq  \K b_1 \cup  \K b_2 \cup \dots \cup  \K b_n $?
Meanwhile, the above problems have been dealt with in \cite{mttmult},
with possible connections with multi-posets \cite{rump} to be further analyzed.

\item  
Characterize those universal-existential
sentences valid in all structures $\mathbf S(X)$,
for all the above languages, and with $X$ varying
either among topological spaces or closure spaces.
  \end{enumerate}  
\end{problemx}  

\begin{remarkx} \labbel{unex}
Of course, if we consider the closure operation 
$ \K $ as part of the language, then 
topological spaces satisfy 
$ \forall xy\  \K (x \vee y) =  \K x \vee  \K y$,
a sentence which is not necessarily valid in closure spaces
(actually a closure space is a topological space if and only if 
it satisfies this sentence, together with a sentence asserting that the
minimal element is a $ \K $-fixed point). 

 Coming back to the language of specialization semilattices,
we observe that Problem \ref{ur}(3)
has surely different solutions, when $X$ varies on
 topological or closure spaces. Indeed, 
the following sentence
\begin{equation}\labbel{ues} 
\forall xyz (
 x \vee y < z   \ \& \  z  \sqsubseteq x \vee y
\Rightarrow 
\exists\ w ( (x < w \ \& \  w \sqsubseteq x) \text{ or }
 (y < w \ \& \  w \sqsubseteq y) ))
  \end{equation}    
is equivalent to a $\forall\exists$
sentence, since the existential quantifier can be moved shortly after 
the universal quantifier; we have kept it in the present position just
to improve readability.
Interpreted in closure spaces, the sentence \eqref{ues} asserts that 
if $x \cup y$ is not closed,
then either $x$ or $y$ is not closed;
in contrapositive,  the union of two closed subsets is closed. 
Hence if $X$ is a closure space such that 
$ z \sqsubseteq \emptyset $ implies 
$ z = \emptyset $, then \eqref{ues}
holds in $\mathbf S(X)$  
if and only if $X$  is a topological space.
Compare also the sentence \eqref{pp}.
 \end{remarkx}   

\begin{remarkx} \labbel{ep}    
Another argument suggesting that 
the language of specialization semilattices might lack expressive power
is that the meet of two elements need not exist in a specialization semilattice,
though meets (=intersections) always exist in the motivating examples.

In this respect, we point out that in a specialization semilattice with
a minimum element $0$ we can express the condition that two elements
have no common nonzero lower bound, which, in the
motivating examples, corresponds to the assertion
that two subsets are disjoint. In many topological definitions 
the relevant point is to have some disjoint (or not disjoint) subsets,
rather than dealing with their exact intersection.
Compare Section \ref{recov} below. 

``Having meet $0$'' is not always preserved by
homomorphisms. On the other hand, preservation holds for the property of
``not having meet $0$''; this does not necessarily mean that 
a meet, say, of $a$  and $b$, exists and is
not $0$: this means that there is some 
nonzero element $c$  such that  $c \leq a $ and $c \leq b$. This is frequently
called the \emph{overlap contact}  relation; see more comments
shortly after Remark \ref{modal}. 
 See \cite{cs,hyp,csb,cp,rsaif} for more details
and for the study of such structures.  
\end{remarkx}
 }

\subsection{Embedding specialization posets} \labbel{sp}

Representation theorems for specialization posets
admit a much simpler proof. 
\arxiv{Notice that not everything from the theory of specialization semilattices
generalizes:
 the notion of additivity from Definition \ref{supcomp}(c) cannot be
even expressed,
in the absence of a join-semilattice operation. }
The following
proposition is a generalization of \cite[Lemma 3.4]{Ja},
to the effect that in the assumptions it is enough
to deal with a closure poset, rather than with
a closure meet semilattice.
Recall that a topological
specialization poset is a specialization poset
of the form $\mathbf P(X)=( \mathcal P(X), {\subseteq} , {\sqsubseteq} )$,
for some topological space $X$.
Recall also the more general Definition \ref{closop}.

\begin{proposition} \labbel{embp}
Every specialization  poset $\mathbf P $ 
can be embedded into    
a topological
specialization poset
in such a way that the embedding
preserves all existing (possibly infinitary) meets 
in $\mathbf P$.  
 \end{proposition}

\begin{proof} 
Suppose that 
$\mathbf P = (P, \leq_P, \sqsubseteq _P)$
is a specialization poset.
For every $a \in P$,
let $ \da a = \{b \in P  \mid b \leq_P a\} $ and
let $S = \mathcal P(P)$.   
The function $\iota$  which assigns 
to $a \in P$ the set   $\da a  \in S$
is  an order-embedding from
 $(P, \leq_P)$ to  
$(S, \subseteq  )$ and $\iota$
 preserves existing (possibly infinitary) meets.
So far, we have just recalled a classical argument, 
e.g. \cite[Ch. I, Theorem 9.9]{Harz}.
 We now define a topology on the set $P$  in such a way that 
$\iota$ is a specialization embedding from 
$\mathbf P$ to $S = \mathcal P(P)$ with the induced specialization. 

If $Y \subseteq P$, set
$\K Y = \{ \, c \in P  \mid \text{there is $d \in Y$
such that }  c \sqsubseteq_P d  \,\}$.
The operation $\K$ is additive by construction;
it is inflationary by \eqref{s4}  
and it is idempotent by \eqref{s3}.
Since $\K \emptyset = \emptyset $,  $\K$ is the closure for some
topology $\tau$  on $P$.  
The induced specialization $ \sqsubseteq _S$ on $S= \mathcal P(P)$ 
turns out to be defined by
$ X \sqsubseteq _S Y$
if, for every $c \in X$,
there is $d \in Y$
such that 
$ c \sqsubseteq_P d $.  

Since $(S, \subseteq, \sqsubseteq _S)$ is topological by 
construction, it remains to show that
$\iota$  is also a $ \sqsubseteq $-embedding.
We need to show that  if $a, b \in P$,
then $ a \sqsubseteq_P b $
if and only if 
$ \da a \sqsubseteq_S \da b $.
Indeed, if $ a \sqsubseteq_P b $
and $c \in \da a$, then 
$c \leq_P a \sqsubseteq_P b$,
thus   $c \sqsubseteq_P b$
by \eqref{s6}. Since $b \in \da b$,
we get  $ \da a \sqsubseteq_S \da b $.
In the other direction, if $ \da a \sqsubseteq_S \da b $,
then, since 
$a \in \da a$, there is 
$d \in \da  b$ such that $ a \sqsubseteq_P d $.
Since  $d \in \da  b$, we have $d \leq_P b$,
hence we get $a \sqsubseteq_P b$ by \eqref{s5}.  
\end{proof} 

The topology $\tau$ on $P$ constructed in the proof of
Proposition \ref{embp} depends only on $ \sqsubseteq $,
not on $\leq$, and is the Alexandrov 
topology associated to $ \sqsubseteq $ \cite[II, 1.8]{J}. 
In particular,  the specialization preorder 
of the topological space $(P, \tau )$ 
is exactly the same specialization relation of the original 
specialization poset $\mathbf P$.  

\arxiv{Notice that, if we just want to prove
that a specialization poset can be embedded
into a specialization semilattice, 
in the proof of Proposition \ref{embp} it is enough to take a much smaller subset
of $\mathcal P(P)$, namely, it is enough to take 
 $S$ equal to the family containing all the  finite unions
of   sets of the form $\da a$, with $a$ varying in $P$
(notice that, by construction, $S$, as  defined here,  is closed under finite unions).

Notice also that we cannot use the argument in the proof of 
Proposition \ref{embp} in order to prove
Theorem \ref{sptop},
since $\iota$ from the proof of   \ref{embp}
is  generally not join preserving.
}

The same proof of Corollary \ref{corun}
provides the following corollary  about the language
of specialization posets.

\begin{corollary} \labbel{corcor}
Suppose that $\sigma$  is a  universal first-order sentence
in the language of specialization posets.
Then  $\sigma$  
holds in all the models of the form
$\mathbf P(X)$ (where $X$ 
can equivalently 
vary among topological spaces, closure spaces, or even
closure posets)
if and only if  $\sigma$
 is a logical consequence
of the theory of specialization posets.
 \end{corollary}

\subsection{The amalgamation property} \labbel{ap}

The proof of Theorem \ref{propcomp} provides  canonical extensions
which can be used to show that the theory of
specialization semilattices has the amalgamation property.
Since the theory is locally finite
(i.e., every finitely generated model is finite), then classical model-theoretical
methods provide the existence of  Fra\"\i ss\'e limits
and of the model-completion. 

\begin{lemma} \labbel{embest}
Suppose that $\mathbf  S$ and $\mathbf A$ 
are specialization semilattices and
$\lambda : \mathbf S \to \mathbf A$ is an
embedding (resp., a homomorphism).
Let $\mathbf  U _{\mathbf  S} $
and $\mathbf  U _{\mathbf  A} $ be
 the principal specialization semilattices constructed in the proof
of Theorem \ref{propcomp}, in the latter case starting the construction with
  $\mathbf A$ in place of $\mathbf S$.
Let $ \kappa _{_ {\mathbf  S}} : \mathbf  S \to \mathbf   U _{\mathbf  S} $
and $ \kappa _{_ {\mathbf  A}} : \mathbf  A \to  \mathbf  U _{\mathbf  A} $
be the embeddings constructed in the proof of Theorem \ref{propcomp}.

Then the embedding (homomorphism)
$\lambda$ can be lifted to an  embedding (homomorphism)
$\lambda^* : \mathbf U _{\mathbf  S} \to   \mathbf  U _{\mathbf  A}$ 
such that $\lambda^* \circ \kappa _{_ {\mathbf  S}} = 
\kappa _{_ {\mathbf  A}} \circ \lambda$.
Moreover, $\lambda^*$ is a $ \K $-homomorphism.    
 \end{lemma}

\begin{proof}
Recall that 
$ \mathbf  U _{\mathbf  S}$ is a quotient of $\mathbf  S \times \mathbf T$,
where $\mathbf T$ is a two-element specialization semilattice with 
domain $ T =  \{ 0, 1 \} $. Moreover, $ \kappa _{_ {\mathbf  S}} $
is obtained by composing the projection with the embedding
$\iota _{_ {\mathbf  S}} : \mathbf  S \to \mathbf  S \times \mathbf T$   defined by 
$\iota _{_ {\mathbf  S}}(s) = (s,0)$. 
Thus if we define $\lambda^\diamond : \mathbf  S \times \mathbf T
\to \mathbf  A \times \mathbf T$ by
$\lambda^\diamond (s,x)= ( \lambda (s), x)$,
we surely have that  $\lambda^\diamond$
is an embedding (homomorphism) and  
$\lambda^\diamond \circ \iota _{_ {\mathbf  S}} = 
\iota _{_ {\mathbf  A}} \circ \lambda$.
If we show that $\lambda^\diamond$
passes to the quotient to an embedding
(homomorphism) $\lambda ^*$  making the following diagram commute
\begin{equation*}
\xymatrix{
{\mathbf  S}  \ar[d]_{ \lambda } \ar[r]^{\iota _{_ {\mathbf  S}}}
   & \mathbf  S \times \mathbf T \ar[d]_{ \lambda ^\diamond}
   \ar[r]^ {\pi  _{_ {\mathbf  S}}} & \mathbf  U _{\mathbf  S} 
   \ar@{-->}[d]_{ \lambda ^*} \\
{\mathbf  A}  \ar[r]^{\iota _{_ {\mathbf  A}}}
   & \mathbf  A \times \mathbf T 
   \ar[r]^ {\pi  _{_ {\mathbf  A}}} & \mathbf  U _{\mathbf  A}
}
\end{equation*}
then  we obtain the desired embedding (homomorphism).
The non trivial case is when
$(s,1) \sim (r,1)$ in $\mathbf  S \times \mathbf T$,
 for $s,r \in S$.
This means $s \sqsubseteq_S r$ and $r \sqsubseteq_S s$,
thus  $ \lambda (s) \sqsubseteq_A \lambda (r)$
 and $ \lambda (r) \sqsubseteq _A \lambda (s)$,
since $\lambda$ is supposed to be a homomorphism.
This means that 
$( \lambda (s),1) \sim ( \lambda (r),1)$ in $\mathbf  A \times \mathbf T$,
hence 
$\lambda^\diamond (s,1)=( \lambda (s),1) \sim ( \lambda (r),1)
=\lambda^\diamond (r,1)$.
This shows that  $\lambda^*$ is well-defined, hence 
a semilattice homomorphism.

We now prove that $\lambda^*$ is a $ \K $-homomorphism.
The proof of Theorem \ref{propcomp} gives
$ \K _{\mathbf  U _{\mathbf  S}} [s,x]= [s,1]$, hence
 $ \lambda^*( \K _{\mathbf  U _{\mathbf  S}} [s,x])
= \lambda^*( [s,1]) =[ \lambda (s), 1]=
 \K _{\mathbf  U _{\mathbf  A}} [ \lambda (s), x] \allowbreak =
\K _{\mathbf  U _{\mathbf  A}}  \lambda^* ([s, x])$.
Since, as we mentioned in Definition \ref{closop}, 
 every semilattice- and $ \K $-homomorphism is a 
$ \sqsubseteq $-\hspace{0 pt}homomorphism,
we get that 
 $\lambda^*$ is a 
$ \sqsubseteq $-homomorphism

Finally, suppose that  $\lambda$ is an embedding.
Then $\lambda^*$ is injective, by the same argument
in the proof that  $\lambda^\diamond$ passes to 
the quotient, going in the reverse direction.
If 
 $\lambda^*([s,x])  \mathrel {\undertilde{\sqsubset}_{U_A}} \lambda^*([r,y])$,
that is,
 $[ \lambda (s),x]  \mathrel {\undertilde{\sqsubset}_{U_A}} [ \lambda (r),y]$,
then  $( \lambda (s),x)\sqsubseteq_{A \times T} ( \lambda (r),y)$,
by the first two lines in the proof of Lemma  \ref{quot} and 
since 
 $ \sim$ satisfies \eqref{seq}. 
From $( \lambda (s),x)\sqsubseteq _{A \times T} ( \lambda (r),y)$ we get
$ \lambda (s) \sqsubseteq _{A} \lambda (r)$ and $x \sqsubseteq_{T} y$, thus
if $\lambda$ is an embedding, then $s \sqsubseteq_{S} r$,
hence  $(s,x)\sqsubseteq _{S \times T} ( r,y)$, thus
$[s,x]  \mathrel {\undertilde{\sqsubset}_{U_S}} [ r,y]$.
This shows that 
$\lambda^*$ is an embedding.
\end{proof}

Recall that a class $\mathcal K$ 
 of models of the same type has the
\emph{joint embedding property}
(the \emph{amalgamation property}) if, 
for all models 
$ \mathbf A, \mathbf B\in \mathcal K$ 
(for all models $\mathbf S, \mathbf A, \mathbf B\in \mathcal K$ 
and
embeddings $ \lambda  \colon \mathbf S \to \mathbf A$
and 
 $ \theta  \colon \mathbf S \to \mathbf B$)
 there are a model
$\mathbf D \in \mathcal K$ and  embeddings
$ \lambda _1 \colon \mathbf A \to \mathbf D$
and 
 $ \theta_1   \colon \mathbf B \to \mathbf D$
(such that 
$  \lambda _1 \circ \lambda
=
\theta _1  \circ \theta $).
\arxiv{If the last conclusion
can be strengthened to get $\lambda_1^\sharr(A) \cap \theta_1^\sharr  (B)
 = (\lambda_1 \circ \lambda)^\sharr(S)$, then
$\mathcal K$ 
 is said to have the
\emph{strong amalgamation property}. }

\begin{theorem} \labbel{closl}
\cite{sapimpap} 
The class of closure semilattices has the 
\arxiv{strong }
 amalgamation property.
 \end{theorem}

 \begin{proof} 
The theorem is proved in \cite{sapimpap}  in a more general context.
We present here a sketch of the proof for
the reader's convenience. Since the class of closure semilattices is
closed under isomorphism, we can
 assume that $ \mathbf  S \subseteq \mathbf A $, 
$ \mathbf  S \subseteq \mathbf B $ and $A \cap B = S$
in the premises
of the amalgamation property,
namely, the embeddings are inclusions.
On $D=A \cup B$ let $\leq_ {\mathbf  D}$  
be the union of $\leq_ {\mathbf  A}$, $\leq_ {\mathbf  B}$,
${\leq_ {\mathbf  A}} \circ {\leq_ {\mathbf  B}}$
and ${\leq_ {\mathbf  B}} \circ {\leq_ {\mathbf  A}}$.
It can be checked \cite{Jo,duerel} that $(D, \leq_ {\mathbf  D})$
\arxiv{strongly }
amalgamates the order-reducts of $\mathbf A$ and $\mathbf  B$
over $\mathbf  S$. Moreover, since $\mathbf S$ 
is a join semilattice, joins of $\mathbf A$ are preserved in
$(D, \leq_ {\mathbf  D})$, and similarly for  joins of $\mathbf B$
\cite[p.\ 205]{Jo}.
Since every poset $\mathbf  D$  can be extended to a  complete lattice in such
a way that existing   joins in $D$  are preserved \cite[Sect.\ 1.10]{Harz},
we can extend $(D, \leq_ {\mathbf  D})$ to a  complete lattice
$\mathbf  E$ which amalgamates the semilattice reducts of
 $\mathbf A$ and $\mathbf  B$
over $\mathbf  S$. So far, we have just repeated some arguments
from   \cite{Jo} used to prove the  amalgamation property 
for lattices; we just note here that the argument works for 
semilattices, as well.
Now define $ \K  $ on $E$ by 
$ \K x = \bigwedge \{ \,  \K _ {\mathbf  A} a \mid 
a \in A, x \leq_ {\mathbf  E}  \K _ {\mathbf  A} a \, \}  \wedge 
\bigwedge \{ \,  \K _ {\mathbf  B} b \mid 
b\in B, x \leq _ {\mathbf  E}  \K _ {\mathbf  B} b \, \}$.
It can be checked that $ \K   $  is a closure operator on $E$
and $ \K $ extends both $ \K _ {\mathbf  A}$ and $ \K _ {\mathbf  B}$.
This is essentially the same argument as in \cite[Lemma 2.3]{MT}.
If we expand $\mathbf  E$ by adding the operation  $ \K $,
we get an amalgamating model for the original closure
semilattices. 
\end{proof}  

If $\mathcal H$  is a class of finitely
 generated models in a countable type, a
\emph{Fra\"\i ss\'e}   limit of $\mathcal H$  is a
 countable universal homogeneous model of age $\mathcal H$.
Recall that the \emph{age}  of some model $\mathbf M$  is the class  of all
finitely generated models that can be embedded in $\mathbf M$.
A model $\mathbf M$ is \emph{homogeneous}  if every isomorphism between
finitely generated submodels of $\mathbf M$  extends to an automorphism of 
the whole  $\mathbf M$.
Classical examples are the following: the Fra\"\i ss\'e limit of 
the class of finite linearly
ordered set is the ordered set  of the rationals.
 The Fra\"\i ss\'e limit
of the class of finite graphs is the  random graph. 
See \cite[Chapter 7]{H}  for details.
Random structures have also been  sometimes
 considered from a philosophical point of view \cite{Mie,Rus}.
 
A first-order theory is \emph{model-complete} if every embedding
between its models is elementary. 
If $T$ and $T^*$ are
first-order theories, 
then $T^*$ is said to be the \emph{model-completion} of $T$  
if  $T$ and $T^*$ have the same universal consequences,
$T$  has  the amalgamation property and
$T^*$  is model-complete. 
A theory is $ \omega$-categorical if all its countably infinite
models are isomorphic.
 See again \cite{H} for details.

\begin{theorem} \labbel{apthm}
The class of specialization semilattices has the 
\arxiv{strong} 
 amalgamation property.

The class of finite 
 specialization semilattices has
 a Fra\"\i ss\'e limit $\mathbf M$.
The first-order theory of $\mathbf M$ is 
$ \omega$-categorical, has quantifier elimination and is  the model-completion
of the theory of specialization semilattices.  
 \end{theorem} 

\begin{proof}
Suppose that $\mathbf S , \mathbf A, \mathbf B$ and
$ \lambda  \colon \mathbf S \to \mathbf A$,
$ \theta  \colon \mathbf S \to \mathbf B$
are as in the assumptions of the  amalgamation property. 
 By Lemma \ref{embest} the embeddings $\lambda$ and
$\theta$ can be lifted to  embeddings 
$\lambda^*$ and
$\theta^*$ such that  in  the 
 diagram below the parts denoted by solid lines commute,
where
 $ \mathbf  U _{\mathbf  S},  \mathbf  U _{\mathbf  A},  \mathbf  U _{\mathbf  B}$ 
are principal specialization semilattices and $\lambda^*$, $\theta^*$
are $ \K $-homomorphisms. This means that if we expand
  $ \mathbf  U _{\mathbf  S},  \mathbf  U _{\mathbf  A},  \mathbf  U _{\mathbf  B}$
by adding the closure operation, then $\lambda^*$ and $\theta^*$
are homomorphisms for the expanded models. 
By Theorem \ref{closl} the expanded models can be  amalgamated 
to some closure semilattice $\mathbf  D$, as pictured at the top
 of the diagram below.
 Since every closure lattice homomorphism
is a homomorphism for the specialization reduct, the
specialization reduct of $\mathbf  D$  amalgamates
$ \mathbf A $ and $  \mathbf B$ over $\mathbf S $
by means of the embeddings  
$ \lambda  _1 \circ \kappa  _{_ {\mathbf  A}}$ and 
 $ \theta _1 \circ \kappa  _{_ {\mathbf  B}}$. 
The above argument proves the  amalgamation property.
\begin{equation*}
\xymatrix{
& \mathbf  D &
\\
   \mathbf  U _{\mathbf  A} \ar@{-->}[ru]^{ \lambda_1}
&&
\mathbf  U _{\mathbf  B}  \ar@{-->}[lu]_{ \theta_1}
\\
\mathbf  A  \ar[u]^{ \kappa  _{_ {\mathbf  A}} } 
&   \mathbf  U _{\mathbf  S} \ar[lu]_{ \lambda^* }  \ar[ru]^{ \theta ^* }
 &
\mathbf  B  \ar[u]_{ \kappa  _{_ {\mathbf  B}}} 
\\
& \mathbf S \ar[lu]^{ \lambda }  \ar[ru]_{ \theta  } \ar[u]_{ \kappa  _{_ {\mathbf  S}}}  &
}\end{equation*}

\arxiv{The \emph{strong} amalgamation property
does not automatically follow from the commutativity of the
above diagram. However, for $a \in A$ and $b \in B$,  
 $ \lambda  _1 ( \kappa  _{_ {\mathbf  A}}(a))=
 \theta _1 (\kappa  _{_ {\mathbf  B}}(b))$
means 
$ \lambda  _1 ( [a,0])=
 \theta _1 ([b,0])$.
By the ``strong'' in Theorem \ref{closl},
$[a,0] = \lambda ^* ([c,d])$ and    
$[b,0] = \theta  ^* ([c,d])$,
for some $[c,d] \in \mathbf  U _{\mathbf  S}$. 
The definitions of $\lambda ^*$ and $ \theta  ^*$ (together
with the definition of the equivalence relation in the proof of 
Theorem  \ref{propcomp}) imply that 
$d=0$ and $a = \lambda (c)$,
with $c \in S$. Thus
 $ \lambda  _1 ( \kappa  _{_ {\mathbf  A}}(a))$
belongs to the image of $S$ under 
$  \lambda  _1 \circ \kappa  _{_ {\mathbf  A}} \circ  \lambda $.
This shows the strong amalgamation property.
 }
 
All the above arguments work also in case $S$ is empty;
in this way we get the joint embedding property.
All the above constructions preserve  finiteness, hence the class of
finite  specialization semilattices
has the amalgamation property and the joint embedding property.
The existence of a Fra\"\i ss\'e limit follows now from
\cite[Theorem 7.1.2]{H}. The last sentence follows from
 \cite[Theorem 7.4.1]{H} and \cite[Fact 2.1(3)]{KS},
since the class of specialization semilattices is uniformly locally finite.
 \end{proof}

\arxiv{Notice that if, in the construction in the proof 
of Theorem \ref{apthm} we start with (the specialization
reducts of) closure semilattices, then the construction
modifies the closure. Compare Remark \ref{prac}(b). }

The theory of specialization posets, too, has the  amalgamation property,
Fra\"\i ss\'e limits etc. The proof is much simpler, essentially,
the proof amounts to employ simultaneously for both
$\leq$ and $ \sqsubseteq $ the argument showing
that the class of posets has the amalgamation property, as sketched
here at the beginning of the proof of Theorem \ref{closl}.
Full details appear in \cite{duerel}.

\section{Recovering some topological notions} \labbel{recov}

Some topological notions are expressible even in the weaker language
of specialization posets. In the present section
 ``specialization poset'' can be always  replaced by
 ``specialization semilattice''. 
We say that an element $c$ of a specialization poset  $\mathbf S$ is
\emph{closed} if, for every $a \in S$,  $ a \sqsubseteq c$ implies $a \leq c$.
In a specialization poset meets need not exist,
but it is immediate to show that if two closed elements do have a meet,
then their meet is closed, as we  prove in
 the next lemma.

\begin{lemma} \labbel{meetclosed}
Suppose that 
$\mathbf S$ is a specialization poset.
If $c$ and  $d$ are closed elements of $\mathbf S$
 and the meet $ c \wedge d$ exists in $\mathbf S$,
then  $ c \wedge d$ is closed.
 \end{lemma}

 \begin{proof} 
If $a \sqsubseteq  c \wedge d$, then  $ a \sqsubseteq c$, by \eqref{s5},
hence  $a \leq c$, since $c$ is closed.
Similarly, $a \leq d$, thus $a \leq c \wedge d$.
This shows that $ c \wedge d$ is closed.
\end{proof}  

In fact, the analogue of Lemma \ref{meetclosed}
holds with the same proof for the meet of any family, possibly infinite, of
closed elements.

A \emph{specialization poset (semilattice) with $0$}
is a specialization poset (semilattice) with a minimum $0$  which   
is closed, namely, $ a \sqsubseteq 0$ implies $a=0$. 
A \emph{homomorphism $\varphi: \mathbf S
\to \mathbf T$  of specialization posets or semilattices
with $0$} is assumed to send $0_{\mathbf S}$
to $0_{\mathbf T}$, and we also require that if 
$a \in S$ and $a \neq 0$, then $\varphi(a) \neq 0$.
Of course, this condition is satisfied in our motivating examples,
where homomorphisms correspond to the image function 
associated to some function.

As we mentioned, meets need not exist in a specialization poset;
however,  in the presence of a  $0$, 
we can express the statement that two elements $a$ and  $b$  have meet 
equal to $0$ (in the
motivating examples, ``the intersection of $a$ and  $b$  is empty'') as follows:
$\forall x (  x \leq a \,\& \, x \leq b \Rightarrow x=0)$.
Similarly, we can express by a possibly infinitary sentence the
assertion that the meet of a, possibly infinite, family of elements is $0$.
This is enough in order to express some 
topological properties, e.g.\ compactness.  \arxiv{(Of course,
here we are temporarily 
renouncing to Requisite (R1) from the introduction!)}

A specialization poset (semilattice) with $0$  is \emph{compact} 
if, for every family $D$ of closed elements with meet $0$,
there is a finite subfamily  $D_F \subseteq D$ with meet $0$.
We can then generalize to specialization posets the
basic topological fact that a closed subset of a compact topological space
is still compact. If $\mathbf S$ is a specialization poset 
and $c \in S$, the \emph{specialization poset  restricted
to $c$} is the submodel of $\mathbf S$ with domain
$S_{\da c} = \{ \, s \in S  \mid s \leq c \,\}$. 
Note that if $\mathbf S$ is a  semilattice,
then $S_{\da c}$ is closed under the semilattice operation,
hence $S_{\da c}$ is the domain for a semilattice.

\begin{lemma} \labbel{cpcl}
Suppose that  $\mathbf S$ is a  specialization poset or semilattice.
  \begin{enumerate}[(1)]    
\item 
Suppose that $d , c \in S$, $d \leq _{\mathbf S}  c$ and  $c$  is closed in $\mathbf S$.
If $d$ is closed in  $\mathbf S_{\da c} $,
then $d$ is closed in  $\mathbf S$.
\item
If $\mathbf S$ has $0$, $\mathbf S$ is  compact 
 and $c \in S$ is closed in $\mathbf S$,
 then $\mathbf S_{\da c} $ is compact.  
   \end{enumerate} 
 \end{lemma} 

\begin{proof} 
(1) If $a \in S$ and $a \sqsubseteq _{ \mathbf S} d$,
then   $a \sqsubseteq _{ \mathbf S} c$, by \eqref{s5}, 
since $d \leq _{\mathbf S} c$. Since $c$ is closed in  $\mathbf S$,
then $a \leq  _{\mathbf S} c$, thus  $ a \in  S_{\da c} $.
Since, by definition, the specialization
$ \sqsubseteq $ on $  \mathbf  S_{\da c} $
is the restriction of $ \sqsubseteq _{ \mathbf  S} $
to $S_{\da c} $, then  from $a \sqsubseteq _{ \mathbf S} d$
we get $a \sqsubseteq d$ in  $  \mathbf  S_{\da c} $,
since $a,d \in S_{\da c}$.
Since $d$ is closed in $  \mathbf  S_{\da c} $,
then $a \leq d$ in   $  \mathbf  S_{\da c} $,
thus $a \leq  _{\mathbf S} d$.

(2) By construction, $S_{\da c}$ is downward closed,
hence a family $D$ of elements of $S_{\da c}$  has meet $0$ in
   $  \mathbf  S_{\da c} $ if and only if 
$D$  has meet $0$ in
   $  \mathbf  S $. The conclusion is now immediate from (1).
\end{proof}    

Under some mild assumptions on  homomorphisms
(assumptions which are all satisfied in the motivating examples),
we  can now generalize the topological theorem 
asserting that the image of a compact topological space
under a continuous function 
is compact.

If $\varphi: \mathbf S \to \mathbf T$ 
is a homomorphism of specialization posets (semilattices) and
$t \in T$,
we say that  \emph{$t$ has maximal preimage}
if in $\mathbf S$ the set $\{ \, s \in  S \mid \varphi (s) \leq _{\mathbf T} t  \,\}$   
has a maximum element.
Notice that if $t$ lies in the image of $S$
under $\varphi$   and $t$ has maximal preimage $s$,
then $\varphi(s)=t$. 
\arxiv{ The next lemma 
generalizes the fact that, for continuous functions between
topological spaces, the preimage of a closed subset is closed. }
  
\begin{lemma} \labbel{preimdiclo} 
Suppose that  $\varphi: \mathbf S \to \mathbf T$ 
is a homomorphism of specialization posets and
$t \in T$ has maximal preimage $s$.

If $t$ is closed in $\mathbf T$, then $s$ is closed in
$\mathbf S$.    
\end{lemma} 

\begin{proof} 
If $a \in S$ and $ a \sqsubseteq _{\mathbf S} s$,
then  $ \varphi (a) \sqsubseteq _{\mathbf T} \varphi  (s) 
\leq  _{\mathbf T} t$, thus $ \varphi (a) \sqsubseteq _{\mathbf T} t$,
by \eqref{s5}, hence  $ \varphi (a) \leq _{\mathbf T} t$,
since $t$ is closed.  
Since $s$ is the maximal preimage of $t$,  
the above inequality means 
$a \leq _{\mathbf S} s$, what we had to show. 
\end{proof}

\begin{theorem} \labbel{imcpt}
Suppose that  $\varphi: \mathbf S \to \mathbf T$ 
is a surjective homomorphism of specialization posets 
(or semilattices) with $0$  and
suppose that 
every closed element of $\mathbf T$ has a maximal preimage in $\mathbf S$.

If $\mathbf S$ is compact, then $\mathbf T$ is compact. 
 \end{theorem} 

\begin{proof}
Suppose that $D$ is a family of closed elements of $\mathbf T$  
and $D$ has meet $0_{\mathbf T}$. For $d \in D$,
let $m(d)$ be the maximal preimage of $d$ and let  
 $E= \{ \, m(d) \mid  d \in D\, \} $.
If $d \in D$ and $c=m(d)$,
then $\varphi(c)=d$, by  the remark  just before the statement of Lemma
 \ref{preimdiclo}
and using the assumption that $\varphi$  is surjective.
Thus the restriction $\varphi_{\restriction E}$  of $\varphi$  to $E$ is a
surjection from $E$ to $D$; it is also injective, since 
$m$ is its inverse.

By the definition of a homomorphism
of specialization posets with $0$, $E$ has meet $0_{\mathbf S}$,
 since otherwise
there is $x \in S$ with  $0_{\mathbf S} <_{\mathbf S} x
 \leq_{\mathbf S}  e$, for every $e \in E$, 
thus   $0_{\mathbf T} < _{\mathbf T} \varphi (x) 
\leq _{\mathbf T} d$, for every $d \in D$, since
$\varphi_{\restriction E}$ is a
surjection from $E$ to $D$. This contradicts the assumption
that $D$ has meet $0$. By Lemma \ref{preimdiclo},
each member of $E$ is a closed element of $\mathbf S$,
thus  $0_{\mathbf S} =  \bigwedge_{\mathbf S} E_F$,
for some finite subset $E_F$ of $E$, since $\mathbf S$ is compact. 

Since we have showed that $\varphi_{\restriction E}$
 is a function from $E$ to $D$,
the set $D_F$ of the images of $E_F$ under $\varphi$  
is a subset of $D$.
 Since $D_F$ is finite, we can conclude the proof if we show that  
 $0_{\mathbf T} =  \bigwedge_{\mathbf T} D_F$.
Suppose not, thus there is $y \in T$ such that 
$0_{\mathbf T}<_{\mathbf T} y
 \leq _{\mathbf T} d$, for every $d \in D_F$.
Since $\varphi$  is surjective, $y = \varphi (x)$,
for some $x \in S$ with $0 _{\mathbf S} < _{\mathbf S} x $,
since we assume that homomorphisms preserve $0$.
 For every $d \in D_F$,
if $c=m(d)$, then 
$c \geq _{\mathbf S} x$, since   $\varphi (x)= y \leq _{\mathbf T} d$
and $c$ is the maximal preimage of $d$.
Since we have showed that $m$  is
the inverse of $\varphi_{\restriction E}$, then $c \in E_F$.
Letting $d$ vary in $D_F$,
we get  $c \geq _{\mathbf S} x$,
for every $c \in E_F$. 
This contradicts $0_{\mathbf S} =  \bigwedge_{\mathbf S} E_F$.
Thus $0_{\mathbf T} =  \bigwedge_{\mathbf T} D_F$ and this shows that 
$\mathbf T$  is compact.
 \end{proof}

\section{Further remarks} \labbel{fr} 

 \begin{remark} \labbel{compl}
If $(S, \leq )$ is a poset, then,
among the relations  $ \sqsubseteq $  making 
 $S$ a specialization poset,
there is obviously the finest 
one, namely,
$ {\sqsubseteq} = {\leq} $,
and there is the coarsest relation, namely 
the universal relation
$ {\sqsubseteq} = \mathcal P(S \times S)$  
such that $a \sqsubseteq b$, for every
$a,b \in S$. 
If in addition $(S, \vee)$ is a
semilattice, then the above coarsest  (finest) relation
makes $S$ a specialization semilattice.

More generally,
given a poset (semilattice) $S$,
the set of all the binary relations making
$S$ a specialization poset (semilattice)
is a complete lattice with maximum and minimum. 
 The maximum and minimum have been described above.
Meet is intersection of relations:
the meet of $( \sqsubseteq _i) _{i \in I} $
is the relation $ \sqsubseteq $
defined by $a \sqsubseteq  b$
if and only if    $a \sqsubseteq _i b$,
for every $i \in I$.
 \end{remark}

\begin{remark} \labbel{modal}
Specialization semilattices can be considered
the algebraization of the fragment of 
(non-normal monotone) modal logics 
consisting of formulas of the type
 $ A \Rightarrow \diamondsuit B$
or $ A \Rightarrow  B$, 
 where $A$ and $B$
are disjunctions of propositional variables
and $\diamondsuit$ is the possibility operator.
Similar fragments have been considered in the literature, 
e.g.\  \cite{Gr,Iv,KKT,MC}.  
Non-normal monotone modal logics 
have recently attracted some interest; see 
\cite{Fr} for a careful history of the subject and
for more references. 

A comparison with \cite{GKWZ,KKT} suggests the problem of
studying semilattices endowed with more than one
specialization. 
 \end{remark}

 There are  more properties 
preserved by image functions, besides the properties we have
considered here. For example, the 
unary relation $A(x)$ expressing the property
that $x$ is an atom (``a singleton''),
and the binary relation $In(x,y)$, \arxiv{frequently called
the \emph{overlap relation}, }
 expressing
the property that  $x$ and $y$ do not have meet  $0$
(``$x$ and $y$ are not disjoint'')
are preserved by image functions.
Proximity notions, too,  are generally
preserved (e.~g., ``$Kx$ and $Ky$ are not disjoint'');
 see \cite{DC} for more details and further references.  
See also Remark \ref{subt}.
\arxiv{Meanwhile, such proximity or ``contact''
relations, as well as some generalizations, have been studied in
\cite{cs,hyp,csb,cp}
from the perspective of the present work.}

\begin{problems} \labbel{prob}
(1) Are there more properties 
preserved by images of continuous functions, in particular, 
properties involving basic topological notions, for example,
 the interior operation?
Is the structure of specialization semilattices and posets
significantly affected if we add corresponding 
relations and axioms? \arxiv{See \cite{rsaif} for merging
various forms of ``contact'' with various forms of
specialization.}

(2) Generalize the results of the present paper
when $\leq$ is only assumed to be
a preorder, or when condition \eqref{s2}
is removed or weakened.
This problem is motivated by e.g.\ 
\cite{BHl, FM, Hol}, respectively,
\cite{Ce}. 
\arxiv{See the next problem  for more details. }
 \end{problems}

\arxiv{

\begin{problem} \labbel{cec}
Study the following notions.

A \emph{\v Cech-poset} is a model
$(P, \leq, {\sqsubseteq} )$ such that 
$(P, \leq)$ is a poset and \eqref{s1},
\eqref{s5}, \eqref{s6} from  Definition \ref{spsem} 
and Remark \ref{conseq} are satisfied. 

A \emph{\v Cech-semilattice} is a model
$(P, \vee, {\sqsubseteq} )$ such that 
$(P, \vee)$ is a semilattice and \eqref{s1}, \eqref{s3}, 
\eqref{s5} and  \eqref{s6}  hold. 

Notice that if $ \K $ is an isotone
 and inflationary---not
 necessarily idempo\-tent---operation
on some poset $(P, \leq)$ (semilattice $(P, \vee)$) then setting, 
as custom by now, $a \sqsubseteq b$ if $a \leq  \K b$,
we get a \v Cech-poset (semilattice).

Note also that properties (equivalent to) \eqref{s5} and \eqref{s6}
appear in many distinct settings, e.~g. \cite[I-1]{GHKLMS},
\cite{KP,CJ}. See \cite{ppio} for more details.  
 \end{problem}

\begin{remark} \labbel{percpn}
(a) A shorter but not constructive
proof of Theorem \ref{propcomp}
(but with the weaker conclusion ``principal'' 
instead of ``principal additive'') 
can be obtained by using the Compactness Theorem.
There is a first-order sentence $\sigma$  asserting that a 
specialization semilattice is principal.
Thus, in order to show that a specialization semilattice $\mathbf S$ 
can be embedded into a principal one
it is enough to prove that $T=Diag(\mathbf S) \cup \{ \sigma  \} $ 
has a model. See \cite[Section 1.4]{H} for the definition and the
basic properties of \emph{diagrams}.  
By compactness, $T$ has a model if and only if
every finite subset $T_F$ of $T$ has a model.
Since every finite subset of $T$ involves only 
a finite number  of (names for) elements of $\mathbf S$,
it is enough to work with the semilattice $\mathbf S'$ generated by 
this finite set of elements.
Since semilattices are locally finite, $\mathbf S'$ is finite,
hence principal, by a remark in Definition  \ref{supcomp}(b), thus a model of $T_F$. 
See  \cite[Remark 4.7(f)]{mttmult}  for further details.

(b) Moreover, as remarked at the beginning of the proof
of Proposition \ref{clemb}, we can assume that $\mathbf S$ has
a minimal element $0$. Since a finite join semilattice   
with a minimum is a lattice, we can add to $T$
the axioms for lattices (in the language of semilattices), hence
 the argument in (a) shows that 
every specialization semilattice  
can be embedded into a principal specialization semilattice which,
as an ordered set, is a
lattice. 
 \end{remark}   

\emph{Conclusions.} 
 While the theory  of specialization
 semilattices presented here 
might prove a bit relevant to foundational studies
about topology, 
it is 
possibly  
too weak
to reproduce an important part of 
topological results.
However,
the theory appears to be interesting for itself,
since it seems to capture 
significant parts of the notions of closure, hull,
generated by\dots, even in the case when
the actual ``closure'' of some set
is too large to be considered  ``admissible''  
in the framework under consideration, or, anyway,
there are reasons suggesting it should  not
necessarily be considered.

The fact that
many examples of this situation appear in many
disparate
unrelated fields of mathematics, with 
applications to other sciences, 
strongly supports the above point of view.
In this sense, our main result asserts that, 
in each of the above situations, we are always
allowed to add ``imaginary''
elements in such a way that we can pretend  to 
be working in an actual topological space.  
Whether or not the above remarks provide
some explanation for the success of 
 topology, the present notions
seem to deserve some study, even 
if they are set (or, possibly, just because they are set)
in an extremely simpler framework.

 }

\arxiv{
\section{Appendix. More examples} \labbel{me}

While we intuitively think of
$a \sqsubseteq b$ as 
``$a$ is contained in the closure of $b$'',
we have provided  examples of 
specialization semilattices which arise in situations
in which no recognizable ``notion of closure'' is present,
e.~g., inclusion
mod finite in
3. in Section \ref{exasub}.  
The next examples also help clarify  
the distinction between homomorphisms and embeddings
in the setting of specialization semilattices.

\begin{example} \labbel{ord}
Let $\alpha$ be an  ordinal,
$\leq$ be the standard order on $\alpha$ 
and define  $\beta \sqsubseteq \gamma $
if and only if there is some natural number $n$
such that $ \beta \leq \gamma +n$.

Then $\mathbf P( \alpha)  = (\alpha, \leq, \sqsubseteq) $
is a specialization poset. Moreover, if we take $\beta \vee \gamma$
to be $\sup \{ \beta, \gamma  \}$, then
 $\mathbf S( \alpha)  = (\alpha, \sup, \sqsubseteq)$ is a specialization semilattice.
In passing, let us notice that \emph{every}
linearly ordered specialization poset
becomes a specialization semilattice, if we consider
the binary
$\sup$ as join.   

(a)    If $\alpha$ is infinite, then
$\mathbf P( \alpha) $ and $\mathbf S( \alpha) $ 
are  not principal.
Recall Definition \ref{supcomp}.
Indeed,  $S_0 = \{ a \in P( \alpha ) \mid  a \sqsubseteq 0 \}$
has no maximum.

(b)
 If $\alpha> \omega $,
then $S_0 $
has a supremum $ \omega$;
however, $ \omega \notin S_0$.  

(c)
 Let $\alpha > \omega $. We now address the following question.

  \begin{enumerate}    
\item[($\diamondsuit$)]  
 Can we give $\alpha$ 
the structure of a principal specialization poset $\mathbf P$ 
in such a way that the identity function is a homomorphism from
$\mathbf P( \alpha) $ to $\mathbf P$?
  \end{enumerate}
Let us look at the possibilities for 
$S_0 = \{ a \in P \mid  a \sqsubseteq 0 \}$, as
evaluated in some hypothetical such $\mathbf P$.
Since we assume that $\mathbf P$ is principal, 
then $S_0$ has a maximum, call it $ \beta $.

If $ \beta $ is not the maximum element of $\alpha$,
that is, if $\beta+1 < \alpha $,   
then $   \beta +1 \sqsubseteq \beta  $
in  $\mathbf P$, since 
$   \beta +1 \sqsubseteq \beta  $
holds in  $\mathbf P(\alpha) $
and we want the inclusion to be a homomorphism.
Since $\beta \in S_0$,
then $\beta \sqsubseteq 0$ in   
 $\mathbf P$, hence 
$   \beta +1 \sqsubseteq 0 $,
by $   \beta +1 \sqsubseteq \beta  $
and \eqref{s2}. This contradicts the assumption that
$\beta$ is the maximum of $S_0$.

Thus the only possibility left is that 
$\alpha= \beta +1$ and $\beta$ is the maximum of
$S_0$, as computed in $\mathbf P$.
Thus in $\mathbf P$ we have
$\gamma \sqsubseteq \delta $,
for every $\gamma, \delta \in \alpha $, by \eqref{s5} and \eqref{s6}.
If $\alpha$ 
is a successor  ordinal, this clearly gives $\mathbf P$ the structure of a 
principal specialization semilattice; moreover, the identity function
is a homomorphism from
$\mathbf P( \alpha) $ to $\mathbf P$, thus 
($\diamondsuit$) has an affirmative answer.
The above arguments show that this is the only way
to accomplish our goal.
Of course, since we have assumed $\alpha > \omega $,
the identity function is not an embedding, since 
$  \omega \not\sqsubseteq 0 $ in  
$\mathbf P( \alpha) $, but 
$  \omega \sqsubseteq 0 $ in  
$\mathbf P $.

(d) On the other hand, by Proposition \ref{embp} 
we can extend   
$ \mathbf P( \alpha) $
 to a principal specialization poset,
and similarly by Theorem \ref{propcomp} we can extend
  $ \mathbf S( \alpha) $
to a principal specialization semilattice.
The above arguments show that in this case  we necessarily should
add new elements.

Following Remark \ref{prac}(a), this can be done by adding
to  $  \alpha $,  for every 
infinite limit ordinal $\gamma \leq \alpha $,
a new element, call it  $\gamma -1$,
where the ordering on the extended set 
  $ S^*( \alpha) $ is defined in the obvious way.
We further set $ \beta \sqsubseteq \gamma {-}1 $
when $\gamma$ is the smallest limit ordinal strictly
larger than $\beta$, together with the further
relations necessary in order to make $ S^*( \alpha) $
a specialization semilattice.

In the notation from the proof of Theorem \ref{propcomp},
 an element $\beta$ 
of  $ S( \alpha) $ is identified with
$[ \beta, 0]$ and $\gamma-1$ corresponds to  
 $[ \beta, 1]$, where, as above, $\gamma$ is the smallest limit ordinal strictly
larger than $\beta$.
 \end{example}    

We now add two scattered remarks 
about some notions
recalled in Section \ref{exasub}.

\begin{remark} \labbel{causs}       
Recall the example of causal spaces from \cite{KP} 
briefly discussed in Section \ref{exasub}, 2b.
The order relations considered in 
\cite{KP} represent \emph{causal precedence}
and \emph{chronological precedence}
 on the points---or ``events''---of a   manifold   modeling
space-time in the  general  theory  of relativity, or, possibly,
in some more abstract generalization.

It might turn out that at the small-scale level 
relativistic events have a composite structure
and some corresponding relations are no more antisymmetric at this
level.
This  might suggest the  shift  from posets to pre-orders.
Of course, under the above interpretation, 
it would be appropriate to consider a theory with two pre-orders,
one finer than the other; in other words, in the definition \ref{spsem}(a)  
of a specialization poset, we should weaken the assumption that $\leq$ 
is an order to a pre-order.
The  antisymmetric 
relations from \cite{KP} would emerge back only at the level of events, after
we take some quotient turning pre-orders into 
orders. 
\end{remark}

\begin{problem} \labbel{2p}
Study the theory of two preorders, one finer than the other.

Is there any significant difference with the theory of 
specialization posets as introduced in Definition \ref{spsem}(a)?
 \end{problem}  

Problem \ref{2p}  is also motivated
by some examples considered
in 2. in  Section \ref{exasub},
in particular, modal frames.

\begin{example} \labbel{meas}
Let $X$ be a set and $\mu $ be a measure defined 
on some subset $S$ of $\mathcal P(X)$.
Then, as already mentioned, $(S, {\subseteq} , {\sqsubseteq} )$
is a specialization poset, where 
$ \sqsubseteq $ is defined by 
\begin{equation}\labbel{mmm}       
\text{$a \sqsubseteq b$ \ \ \  if  \ \ \ 
$\mu (a) \leq \mu (b)$,}
\end{equation}
for $a,b \in S$.

In general, the above definition 
does not furnish a specialization semilattice, 
since  the properties of a measure are incompatible
with \eqref{s3}. 
There is a notable exception:
if $\mu $ is a two-valued measure,
then  
 $(S, \cup , {\sqsubseteq} )$
is indeed a specialization semilattice. 

There is a vast literature on binary relations
representable by the formula \eqref{mmm}, in various contexts,
mostly related with foundational issues about probability and with possible
economical applications. 
See e.g.\ \cite{L} and further references there. 
According to \cite{A}, the whole line of research  dates
to \cite{F}.  
 \end{example}

We now
present a 
few more
counterexamples.

A poset in which some binary joins fail to exist
 cannot be endowed with the structure of 
a specialization semilattice. 
In Example \ref{diamond} 
we have showed  that \eqref{s3} does not follow from the other
assumptions.
We now show that
the assumption that $ \sim$ 
satisfies the condition \eqref{seq} 
in Lemma \ref{quot} is necessary.

 \begin{example} \labbel{exquot} 
Let
$\mathbf S =(S, \vee, {\sqsubseteq} )$, where $S = \{ 0, 1, 2, 3 \}$,
$\vee = \sup$ and  the only nontrivially
$ \sqsubseteq $-related pairs are given by 
$ 1 \sqsubseteq 0 $ and $ 3 \sqsubseteq 2$.
Of course, we also assume 
$ m \sqsubseteq n $, for $m \leq n \leq 3$.  
Then $\mathbf S$ is a specialization semilattice.

If $ \sim$ is the equivalence relation whose classes
are $\{ 0\}, \{ 1,2 \},  \{ 3 \} $, then,
in the notations from Lemma \ref{quot},  
the standard way to define the quotient structure
$\undertilde{\mathbf  S}$
is to set
$ \undertilde{a}  \mathrel {\undertilde{\sqsubset}} \undertilde{b}   $
if and only if  $a_1 \sqsubseteq b_1$,
for \emph{some}
$a_1, b_1$ with  $a_1 \sim a$ and $b_1 \sim b$.
This is required if we want  the
projection to be  a homomorphism.
However $\mathrel {\undertilde{\sqsubset}}$
is not transitive, since  
 $ \undertilde{3}  \mathrel {\undertilde{\sqsubset}} \undertilde{2}
= \undertilde{1}  \mathrel {\undertilde{\sqsubset}} \undertilde{0} $,
but 
  $ \undertilde{3}  \mathrel {\undertilde{\sqsubset}} \undertilde{0}$
does not hold.
Hence the assumption that $ \sim$ 
satisfies the condition \eqref{seq} 
in Lemma \ref{quot} is necessary.
More generally, we have showed that a quotient of 
a specialization semilattice is not necessarily a specialization semilattice.

The point is that, under the assumptions
in Lemma \ref{quot}, in the quotient we are always allowed to choose
the same representative for the ``middle'' element
$\undertilde{b}$ in the implication \eqref{s2}, but this is not
always true
in the general case.   
\end{example}

\begin{example} \labbel{star}
Let $S= \{ a, b, c, 1 \} $
with the partial order $\leq$  
given by 
$a < c < 1$ and $ b < c < 1$.
 Let the only nontrivial
$ \sqsubseteq $-relation be  
$ 1 \sqsubseteq c$.

Then 
$\mathbf S=(S, \vee, {\sqsubseteq} )$
is a principal specialization semilattice,
with $ \K a=a$, $ \K b=b$
and $ \K c= \K 1=1$.
See Definition \ref{supcomp}.    
On the other hand, 
 $\mathbf S$ is not additive, since
$ \K (a \vee b) =  \K c=1> c = a \vee b =  \K a \vee  \K b$.

\begin{center}
\begin{tikzpicture}[-,>=stealth',auto,node distance=3cm,
thick,main node/.style={circle,scale=0.3,draw},
mnode/.style={scale=0.3}]
\node[main node,label={west:$c$}] (c){};
\node[main node,label={north:$\ \ \ \ \ \ \ 1= \K c$}] (1) [above of=c]{};
\node[main node,label={west:$ \K a=a$}] (a) [below left of=c]{};
\node[main node,label={east:$b= \K b$}] (b) [below right of=c]{};
\node[main node,label={west:$a_1$}] (a1) [above left of=a]{};
\draw (a) -- (c) -- (1);
\draw  (b) -- (c);
\draw (a) [dotted] -- (a1) -- (1);
\end{tikzpicture}
\end{center}

If we add a new element $a_1$ to $S$,
prescribing $a < a_1 < 1$
and $  a_1 \sqsubseteq a$, then the resulting structure
$\mathbf S_1$ is a   principal additive
specialization semilattice, since then
 $ \K _1a= a_1 $. 
Moreover, the inclusion is an embedding
of specialization semilattices
(caution! not an embedding, not even a homomorphism, with respect to
the operations  $ \K $, $ \K _1$).

However, we could perform a symmetric construction
by adding some element 
$b_1 >  b$. 
This shows that 
we do not necessarily have a smallest 
extension satisfying  the conclusions of Theorem \ref{propcomp}.
 \end{example} 

Possibly, a more complete model-theoretical approach
to topology can be obtained by forgetting about the requisite (R3)
from the introduction
and considering at the same time properties that 
are preserved covariantly and contravariantly.

\begin{problem} \labbel{coco}
Study the following category $\mathfrak C$. 

Objects of $\mathfrak C$ are closure algebras.

If $\mathbf A$ and $\mathbf  B$ are closure algebras, a 
$\mathfrak C$-morphism from 
$\mathbf A$ to $\mathbf  B$
is a pair of functions $f_*:A \to B$ 
and $f^*:B \to A$ such that 
  \begin{enumerate}[(i)]  
   \item 
$f^*(f_*(a)) \geq a $ and  $f_*(f^*(b)) \leq b $,
for every $a \in A$ and $b \in B$;
\item
$f^*$ is a Boolean homomorphism such that 
$f^*(c)$ is closed in $\mathbf A$, for every 
closed $c$ in $\mathbf  B$;
\item
$f_*$ is order preserving, hence a join-semilattice homomorphism
\cite[Proposition 3.26]{E}. 
  \end{enumerate} 

Note that it follows from (i)-(iii) that 
$f_*$ is continuous, hence a specialization homomorphism,
by Corollary \ref{corr}.  
To show that $f_*$ is continuous, compute
$f_*(Ka) \leq f_*(K f^*( f_*(a))) \leq f_*(f^*(K  f_*(a)))
\leq K  f_*(a)$. 

The above definition can be modified in several ways.
For example,
we can consider Boolean algebras
with a (not necessarily additive) closure operator.
We can actually give a similar definition for just closure
posets, semilattices or lattices in place of closure algebras.
In such cases, it is appropriate to add also 
an interior operator satisfying suitable properties,
since if we lack complementation, closure and interior
are not interdefinable anymore.
 \end{problem}

}

\paragraph{Acknowledgements.}
We express our gratitude to anonymous referees
for a careful reading of the manuscript, for
many useful comments and suggestions, in particular,
suggestions for further references.
The advices of the referees have led to significant improvements
of the manuscript, in particular making it more concise and historically
 more  complete.

\smallskip 

This work has been performed under the auspices of G.N.S.A.G.A.
and has been 
partially supported by PRIN 2012 ``Logica, Modelli e Insiemi''.
The author acknowledges the MIUR 
Department 
Project awarded to the
Department of Mathematics, University of  Tor Vergata, CUP
E83C18000100006.

\bigskip 

\arxiv{
While we have pursued a quite extensive bibliographic search,
and the anonymous referees have kindly provided additional useful references,
the arguments treated or touched upon in the present note are so vast 
that the following list necessarily turns out to be incomplete.
Additional references might be found in the quoted works.

 }

\arxiv{
\renewcommand\refname{Additional References}

 }

{\AuthorAdressEmail{Paolo Lipparini}
{Dipartimento di Matematica
Universit\`a di Roma ``Tor Vergata'' 
\\ 
Viale della  Ricerca Scientifica Specializzata
\\I-00133 ROME ITALY}
{lipparin@axp.mat.uniroma2.it}
}


\begin{thebibliography}{10}    

\bibitem{Ar} Artemov, S., \emph{Modal logic in mathematics},
 in Blackburn, P., van Benthem, J.,  Wolter, F. (eds.), \emph{Handbook of modal logic},
 Stud. Log. Pract. Reason. \textbf{3}, 927--969,
Elsevier B. V., Amsterdam (2007). 

\bibitem{BF}  Barwise, J., Feferman, S. (eds.),
\emph{Model-theoretic logics},
{Perspectives in Mathematical Logic},
{Springer-Verlag, New York}
(1985).

\bibitem{vBB} 
van Benthem, J., Bezhanishvili, G.,
\emph{Modal logics of space}, in
Aiello, M.,  Pratt-Hartmann, I., van Benthem, J. (eds.),
\emph{Handbook of spatial logics},
217--298, Springer, Dordrecht (2007).

\bibitem{B} Bergman, C.,
\emph{Universal algebra. Fundamentals and selected topics},
{Pure and Applied Mathematics (Boca Raton)}
\textbf{301},
CRC Press, Boca Raton, FL
(2012).

\bibitem{BBI}  Bezhanishvili, G.,  Bezhanishvili, N., Iemhoff, R.,
\emph{Stable canonical rules},
{J. Symb. Log.}
\textbf{81}, 284--315 (2016).


\bibitem{BBLM}
 Bezhanishvili, G.,  Bezhanishvili, N., Lucero-Bryan,
  J., van Mill, J.,
\emph{The {M}c{K}insey-{T}arski theorem for locally compact ordered
 spaces},
{Bull. Symb. Log.}
\textbf{27}, 187--211  (2021).

\bibitem{BHl} Bezhanishvili, G., Holliday, W. H.,
 \emph{Locales, nuclei, and {D}ragalin frames}, in
Beklemishev, L., Demri, S., M\'{a}t\'{e}, A.\ (eds.),
\emph{Advances in modal logic} \textbf{11},
 177--196,
Coll. Publ., London
(2016).


\bibitem{BH}
Bezhanishvili, G., Holliday, W. H.,
 \emph{A semantic hierarchy for intuitionistic logic},
{Indag. Math. (N.S.)}
\textbf{30}, 403--469 (2019).  
 

\bibitem{BMM}
Bezhanishvili, G.,  Mines, R., Morandi, P. J.,
\emph{Topo-canonical completions of closure algebras and {H}eyting
              algebras},
{Algebra Universalis}
\textbf{58},
1--34 (2008).
  

\bibitem{Bl} Blass, A.,
 \emph{Combinatorial cardinal characteristics of the continuum},
 in Foreman, M., Kanamori, A. (eds.),
\emph{Handbook of set theory}, 395--489,
 Springer, Dordrecht 
  (2010).


\bibitem{Bly}  Blyth, T. S.,
\emph{Lattices and Ordered
Algebraic Structures},
Springer-Verlag  (2005).
 
\bibitem{CLM}  Caspard, N., Leclerc, B., Monjardet, B.,
\emph{Finite ordered sets. Concepts, results and uses},
{Encyclopedia of Mathematics and its Applications}
\textbf{144}, Cambridge University Press, Cambridge (2012).


\bibitem{CM} Caspard, N., Monjardet, B.,
 \emph{The lattices of closure systems, closure operators, and
              implicational systems on a finite set: a survey},
Discrete Appl. Math.
 \textbf{127},
241--269  (2003).
 

\bibitem{Ce} \v{C}ech, E.,
\emph{Topological spaces}
(Revised edition by Frol\'{\i}k, Z., and Kat\v{e}tov, M.,
 Scientific editor,  Pt\'{a}k, V.,
              Editor of the English translation,  Junge, C. O.),
Publishing House of the Czechoslovak Academy of Sciences,
              Prague; Interscience Publishers John Wiley \& Sons, London-New
              York-Sydney
 (1966).

 
\bibitem{Cel} Celani, S.,
\emph{Quasi-modal algebras},
{Math. Bohem.} \textbf{126},
721--736 (2001).
 


\bibitem{CKc} 
Chang, C., Keisler, H. J.,
 \emph{Continuous model theory},
Annals of Mathematics Studies \textbf{58},
{Princeton Univ. Press, Princeton, N.J.}
(1966).



\bibitem{CPZ}  Conradie, W., Palmigiano, A.,  Zhao,
 Z.,
\emph{Sahlqvist via translation},
{Log. Methods Comput. Sci.}  \textbf{15}, Paper No. 15, 35
(2019).


\bibitem{C} 
 Cs\'{a}sz\'{a}r,  \'{A}.,
\emph{Generalized topology, generalized continuity}, Acta Math. Hungar.
\textbf{96}, 351--357 (2002). 

\bibitem{DIM} 
Delzell, C. N.,
Ighedo, O.,
Madden, J. J.,
\emph{Conjunctive join-semilattices},
{Algebra Universalis}
\textbf{82},
Paper No. 51, 30 (2021).

\bibitem{DC} Di Concilio, A.,
\emph{Proximity: a powerful tool in extension theory, function
 spaces, hyperspaces, {B}oolean algebras and point-free
              geometry},
in Mynard, F.,  Pearl E. (eds), 
\emph{Beyond topology},
  Contemp. Math.
\textbf{486},  89--114,
Amer. Math. Soc., Providence, RI
(2009).


 

\bibitem{DH} Dow, A., Hart, K. P.,
\emph{The
\v Cech-Stone  Compactifications of $ \mathbb N$ 
and $\mathbb R$}, Chapter d-18  
in Hart, K. P.,  Nagata, J.,  Vaughan, J. E. (eds.),
\emph{Encyclopedia of general topology},
Elsevier Science Publishers, B.V., Amsterdam
(2004).

\bibitem{Ed} Edgar, G. A.,
 \emph{The class of topological spaces is equationally definable},
{Algebra Universalis}
\textbf{3},
139--146 (1973).
 

\bibitem{En} Engelking, R.,
\emph{General topology},
Mathematical Monographs
\textbf{60},
PWN---Polish Scientific Publishers, Warsaw
(1977).


\bibitem{E} Ern\'{e}, M.,
 \emph{Closure},
in Mynard, F.,  Pearl E. (eds), 
\emph{Beyond topology},
  Contemp. Math.
\textbf{486}, 163--238,
Amer. Math. Soc., Providence, RI
(2009).


\bibitem{Es} Esakia, L.,
\emph{Heyting algebras},
{Trends in Logic---Studia Logica Library}
\textbf{50},
Translated from the Russian edition (1985) by Anton
 Evseev,
{Springer, Cham}
(2019).

\bibitem{FM}  Fairtlough, M., Mendler, M.,
 \emph{Propositional lax logic},
{Inform. and Comput.}
\textbf{137}, 1--33 (1997).

\bibitem{Fa}
Farah, I.,
 \emph{Luzin gaps},
{Trans. Amer. Math. Soc.}
\textbf{356}, 2197--2239 (2004).


\bibitem{Fr} Frittella, S., 
\emph{Monotone Modal Logic \& Friends}, Ph.D. thesis, Universite d'Aix-Marseille, 2014.  

\bibitem{GKWZ} 
Gabbay, D. M.,  Kurucz, A., Wolter, F., Zakharyaschev,
              M.,
 \emph{Many-dimensional modal logics: theory and applications},
{Studies in Logic and the Foundations of Mathematics}
\textbf{148}, North-Holland Publishing Co., Amsterdam
(2003). 


\bibitem{GM}
 Gabbay, D. M., Maksimova, L., 
\emph{Interpolation and definability. Modal and intuitionistic logics},
Oxford Logic Guides
\textbf{46}, The Clarendon Press, Oxford University Press, Oxford
(2005).


\bibitem{GT} Galatos, N., Tsinakis, C.,
\emph{Equivalence of consequence relations: an order-theoretic and
  categorical perspective},
{J. Symbolic Logic}
\textbf{74},
780--810 (2009).



\bibitem{GHKLMS} 
Gierz, G., Hofmann, K. H., Keimel, K., Lawson, J. D.,
 Mislove, M.,  Scott, D. S.,
 \emph{Continuous lattices and domains},
Encyclopedia of Mathematics and its Applications \textbf{93},
Cambridge University Press, Cambridge
(2003)

\bibitem{GH} 
Goldblatt, R., Hodkinson, I.,
\emph{Spatial logic of tangled closure operators and modal
  mu-calculus},
{Ann. Pure Appl. Logic}
\textbf{168},
1032--1090 (2017).


\bibitem{G} Gr\"atzer, G., \emph{Lattice theory: foundation},
 Birkh\"auser/Springer Basel AG, Basel, 2011. 


\bibitem{Gr}
de Groot, J.,
\emph{Positive monotone modal logic},
{Studia Logica}
\textbf{109},  {829--857},
(2021).

\bibitem{Ha}  Hartshorne, R.,
\emph{Algebraic geometry},
Graduate Texts in Mathematics \textbf{52},
Springer-Verlag, New York-Heidelberg
(1977).

\bibitem
{Harz}
Harzheim, E., \emph{Ordered sets},  
Advances in Mathematics \textbf{7},  New York (2005).

\bibitem{H} Hodges, W.,
\emph{Model theory},
Encyclopedia of Mathematics and its Applications
\textbf{42},
Cambridge University Press, Cambridge
(1993).

\bibitem{Hol} Holliday, W. H.,
\emph{Three roads to complete lattices: orders, compatibility,
 polarity},
Algebra Universalis \textbf{82},
Paper No. 26, 14 (2021).

\bibitem{Ie} Iemhoff, R.,
\emph{Consequence relations and admissible rules},
 {J. Philos. Logic}
 \textbf{45}, 327--348  (2016).

\bibitem{Iv} Ivanova, T.,
\emph{Contact join-semilattices},
{Studia Logica}
\textbf{110}, 1219--1241 (2022).
 
\bibitem{Ja} Jackson, M.,
\emph{Semilattices with closure},
Algebra Universalis
\textbf{52},
1--37 (2004).



\bibitem{Jan} Jansana,  R., 
\emph{Propositional Consequence Relations and Algebraic Logic},
 The Stanford Encyclopedia of Philosophy (Spring 2011 Edition),
 Edward N. Zalta (ed.), 
https://plato.stanford.edu/archives/spr2011/entries/consequence-algebraic/  



\bibitem{Je} 
Je\v{r}\'{a}bek, E.,
\emph{Canonical rules},
  {J. Symbolic Logic}
  \textbf{74}, 1171--1205  (2009).


\bibitem{J} Johnstone, P. T.,
\emph{Stone spaces},
Cambridge Studies in Advanced Mathematics
\textbf{3},
Cambridge University Press, Cambridge
(1982).
 

\bibitem{Jo}
 J{\'o}nsson, B.,
\emph{Universal relational systems},
Math. Scand. \textbf{4}, 193--208  (1956).   

\bibitem{JT} 
J\'{o}nsson, B., Tarski, A.,
\emph{Boolean algebras with operators. {I}},
Amer. J. Math. \textbf{73},
891--939 (1951).

\bibitem{KS} Kaplan, I., Simon, P.,
\emph{Automorphism groups of finite topological rank},
{Trans. Amer. Math. Soc.}
\textbf{372},
2011--2043 (2019).


\bibitem{KKT} Kikot, S., Kurucz, A., Tanaka, Y.,
 Wolter, F., Zakharyaschev, M.,
\emph{Kripke completeness of strictly positive modal logics over
 meet-semilattices with operators},
{J. Symb. Log.}
\textbf{84}, 533--588
(2019).

\bibitem{KP} Kronheimer, E. H., Penrose, R.,
\emph{On the structure of causal spaces},
Proc. Cambridge Philos. Soc.
\textbf{63},
481--501 (1967).

\bibitem{K} 
Kuratowski, C.,
\emph{Topologie. {I}. {E}spaces {M}\'{e}trisables, {E}spaces {C}omplets},
Monografie Matematyczne \textbf{20},
2d ed,
Warszawa-Wroc\l aw
(1948).

\bibitem{La} Lando, T.,
\emph{First order {$S4$} and its measure-theoretic semantics},
{Ann. Pure Appl. Logic}
\textbf{166}, 187--218  (2015).


\bibitem{LV} L\'{a}vi\v{c}ka, T.,  Verner, J. L.,
\emph{Completely separable {MAD} families and the modal logic of
              {$\beta\omega$}},
J. Symb. Log. \textbf{87}, 498--507
(2022).


\bibitem{L} Lehrer,  E.,
   \emph{On a representation of a relation by a measure},
   J. Math. Econom. \textbf{20}, 107--118  (1991).


 \bibitem{ecca} Lipparini, P., 
 \emph{Existentially complete closure algebras},
Boll. Un. Mat. Ital. D (6)
\textbf{1}, 13--19 (1982).

\bibitem{mttlib}  Lipparini, P., 
\emph{Universal extensions of specialization semilattices}, 
Categ. Gen. Algebr. Struct. Appl. \textbf{17}, 
101--116  (2022).

\bibitem{mttna}
  Lipparini, P., 
\emph{Universal specialization semilattices}, 
{Quaestiones Mathematicae} 46, 2163-2176, 2022.  

\bibitem{sapimpap} Lipparini, P.,  \emph{Preservation of superamalgamation by expansions}, arXiv:2203.10570v3,  1--20 (2022).

\bibitem{duerel} Lipparini, P.,  \emph{Comparable binary relations and the amalgamation property}, arXiv:2301.12482,  1--13 (2022).

\bibitem{MC} Mench\'{o}n, P.,  Celani, S.,
\emph{Monotonic modal logics with a conjunction},
{Arch. Math. Logic}
\textbf{60}, 857--877 (2021).
 

\bibitem{MT}  McKinsey, J. C. C., Tarski, A., \emph{The algebra of
   topology}, Ann. of Math.  \textbf{45}, 141--191 (1944).

\bibitem{MT2}  McKinsey, J. C. C., Tarski, A., 
\emph{Some theorems about the sentential calculi of {L}ewis and
{H}eyting},
J. Symbolic Logic \textbf{13}, 1--15
 (1948).

\bibitem{Mie} Mielants, W.,
\emph{Can the best of all possible worlds be a random structure?},
Log. Anal., Nouv. S{\'e}r. \textbf{36},
61--73
(1993).


\bibitem{MN} Montalb\'{a}n, A., Nies, A.,
\emph{Borel structures: a brief survey}, in
Greenberg, N., Hamkins, J. D.,
 Hirschfeldt, D., Miller, R. (eds.), 
\emph{Effective mathematics of the uncountable},
Lect. Notes Log.
\textbf{41}, 124--134,
Assoc. Symbol. Logic, La Jolla, CA
(2013).

\bibitem{MP} 
Mynard, F.,  Pearl E. (eds.), 
\emph{Beyond topology},
  Contemp. Math.
\textbf{486},
Amer. Math. Soc., Providence, RI
(2009).

\bibitem{Nag}
Nagata, J.,
\emph{Modern general topology},
 North-Holland Mathematical Library
\textbf{33},
 Second revised edition,
 North-Holland Publishing Co., Amsterdam
(1985). 

\bibitem{Now} Nowak, M.,
\emph{A syntactic approach to closure operation},
{Bull. Sect. Logic Univ. \L \'{o}d\'{z}}
\textbf{46},
219--232 (2017).
 

\bibitem{No} Nowak, M.,
\emph{Disjunctive and conjunctive multiple-conclusion consequence
        relations},
{Studia Logica}
\textbf{108}, 1125--1143 (2020).


\bibitem{Pp}  Panangaden, P.,
\emph{Causality in physics and computation},
{Theoret. Comput. Sci.},
 \textbf{546}, 10--16 (2014).



\bibitem{PW} Peters, J. F., Wasilewski, P.,
\emph{Tolerance spaces: Origins, theoretical aspects and applications},
Inform. Sci. \textbf{195}, 211--225  (2012).


\bibitem{PP} Picado, J.,  Pultr, A.,
\emph{Frames and locales. Topology without points},
 Frontiers in Mathematics,
 Birkh\"{a}user/Springer Basel AG, Basel,
(2012). 

\bibitem{PP2} Picado, J.,  Pultr, A.,
\emph{Separation in point-free topology},
 Birkh\"{a}user/Springer, Cham,
(2021).

\bibitem{R} Ranzato, F.,
\emph{Closures on {CPO}s form complete lattices},
Inform. and Comput. \textbf{152},
236--249 (1999).


\bibitem{RS}
Rasiowa, H., Sikorski, R.,
\emph{The mathematics of metamathematics},
  Monografie Matematyczne \textbf{41},
Pa\'{n}stwowe Wydawnictwo Naukowe, Warsaw
(1963).
 

\bibitem{Rus} Russell, J. S., 
\emph{On the Probability of Plenitude},
Journal of Philosophy \textbf{117}, 
267--292 (2020)


\bibitem{SG} Sambin, G., Gebellato, S.,
\emph{A preview of the basic picture: a new perspective on formal
              topology}, in
Altenkirch, T.,  Naraschewski, W., Reus,
B. (eds.),
  \emph{Types for proofs and programs ({I}rsee, 1998)},
 {Lecture Notes in Comput. Sci.}
 \textbf{1657},
 194--207,
 Springer, Berlin
 (1999).
 

\bibitem{San}  Santocanale, L.,
\emph{A duality for finite lattices}, HAL 00432113, 1--37 (2009).

\bibitem{Sc} Scowcroft, P.,
\emph{Existentially closed closure algebras},
Notre Dame J. Form. Log. \textbf{61},
623--661 (2020).


\bibitem{Se} Servi, M.,
\emph{Algebre di {F}r\'{e}chet: {U}na classe di algebre booleane con
              operatore},
{Rend. Circ. Mat. Palermo (2)}
\textbf{14},
335--366 (1965).
 

 \bibitem{Si} Sikorski, R.,
\emph{Closure algebras},
Fund. Math. \textbf{36},
165--206
(1949).

\bibitem{Sla}   
\v{S}lapal, J.,
\emph{On categories of ordered sets with a closure operator},
{Publ. Math. Debrecen},
\textbf{78}, 61--69 (2011).

\bibitem{Sto} 
Stone, M. H.,
\emph{The theory of representations for {B}oolean algebras},
Trans. Amer. Math. Soc. \textbf{40},
37--111 (1936).


\bibitem{T} Tarski, A., 
\emph{Logic, Semantics, Metamathematics. Papers from 1923 to 1938},
 J. Corcoran (ed.), Hackett Publishing Co., Indianapolis, IN
 (1983).

\bibitem{U} 
Urquhart, A.,
\emph{A topological representation theory for lattices},
{Algebra Universalis} \textbf{8},
45--58 (1978).


\bibitem{V} Vickers, S.,
\emph{Topology via logic},
Cambridge Tracts in Theoretical Computer Science
\textbf{5},
Cambridge University Press, Cambridge
(1989).


\bibitem{W} W\'{o}jcicki, R.,
\emph{Theory of logical calculi. Basic theory of consequence operations},
Synthese Library
\textbf{199},
Kluwer Academic Publishers Group, Dordrecht
(1988).


\bibitem{Z} Ziegler,  M., \emph{Topological model theory},
in Barwise, J., Feferman, S. (eds.),
\emph{Model-theoretic logics},
Perspectives in Mathematical Logic,
  557--577,
Springer-Verlag, New York
 (1985). 

\end{thebibliography}

\begin{thebibliography}{CK}    

\bibitem[AHS]{joy}
Ad\'{a}mek, J.,  Herrlich, H.,  Strecker, G. E.,
 \emph{Abstract and concrete categories. The joy of cats.},
 {Pure and Applied Mathematics (New York)},
John Wiley \& Sons, Inc., New York
 (1990).

\bibitem[A]{A}  Armstrong,  T. E.,  \emph{Review of \cite{L}},
MR1068224  (1991).  

\bibitem[CJ]{CJ} Celani, S., Jansana, R.,
  \emph{Bounded distributive lattices with two subordinations},
 in 
Mojtahedi, M., Rahman, S., Zarepour, M. S.,
(eds.), \emph{Mathematics, logic, and their philosophies---essays in honour
              of {M}ohammad {A}rdeshir}, 
{Log. Epistemol. Unity Sci.}
\textbf{49},  217--252,
Springer, Cham (2021). 


\bibitem[CK]{CK} 
Chang, C. C., Keisler, H. J.,
 \emph{Model theory},
Studies in Logic and the Foundations of Mathematics \textbf{73},
{North-Holland Publishing Co., Amsterdam-London; American
 Elsevier Publishing Co., Inc., New York} (1973),
third expanded edition (1990).

\bibitem[dF]{F} de Finetti, B.,
\emph{La pr\'evision: ses lois logiques, ses sources subjectives},
Annales de l'I. H. P. \textbf{7},  1--68 (1937).  


\bibitem[F]{fiore}  Fiore, C.,
\emph{Reading Conclusions Conjunctively}, 
Journal of Philosophical Logic \textbf{53},1641--1672 (2024). 

\bibitem[L1]{cs} Lipparini, P., \emph{Contact semilattices},
{Logic Journal of the IGPL} \textbf{32},
815--826 (2023).

\bibitem[L2]{mttmult} Lipparini, P., 
\emph{Multi--argument specialization semilattices}, Mathematica Bohemica,
1--23 (2025).

\bibitem[L3]{hyp}  
Lipparini, P., \emph{Hypercontact semilattices}, Journal of Applied Non-Classical Logics, 1--26  (2025).  https://doi.org/10.1080/11663081.2025.2452738 


\bibitem[L4]{csb} Lipparini, P.,
 \emph{Contact join-semilattices are not finitely axiomatizable},
Logic Journal of the IGPL \textbf{33}, 1--12  (2025). 


\bibitem[L5]{cp} Lipparini, P., \emph{Contact posets}, 
arXiv:2303.06259v2, 1--8  (2023).

\bibitem[L6]{ppio} Lipparini, P.,
\emph{Pairs of partial orders and the amalgamation property},
arXiv:2303.06205, 1--18 (2023). 

\bibitem[L7]{mttcong} Lipparini, P.,
\emph{Semilattices with a congruence}   
arXiv:2407.08446, 1--8 (2024).

 \bibitem[L8]{rsaif} Lipparini, P.,
\emph{Relational structures associated to topological spaces
 and preserved by image functions}, Technical Report,
1--28,  University of Rome Tor Vergata (2025).
https://art.torvergata.it/handle/2108/415551 


\bibitem[Pa]{Pa} Pasini, A.,
\emph{On the {$1\deg $}-order translations of the theory of the
  geometric closure structures},
 {Boll. Un. Mat. Ital. B (5)}
\textbf{18},
 217--230 (1981). 


\bibitem[PN]{PN}
Peters, J., Naimpally, S.,
\emph{Applications of near sets},
Notices Amer. Math. Soc.
\textbf{59},
536--542 (2012).

\bibitem[Ru]{rump} Rump,  W.,
\emph{Multi-posets in algebraic logic, group theory, and non-commutative topology},
Ann. Pure Appl. Logic \textbf{167}, 1139--1160  (2016). 

\end{thebibliography}
\end{document}